\documentclass[12pt, a4paper]{article}
\usepackage{epsfig,graphicx}
\usepackage[utf8]{inputenc}
\usepackage{amsmath}
\usepackage{amsfonts}
\usepackage{amssymb}

\usepackage{color}
\usepackage{hyperref}

\usepackage[dvipsnames]{xcolor}

\usepackage{amsmath}
\usepackage{tikz}
\usepackage{mathdots}

\usepackage{cancel}
\usepackage{color}

\usepackage{array}
\usepackage{multirow}
\usepackage{amssymb}
\usepackage{gensymb}
\usepackage{tabularx}

\usepackage{booktabs}
\usetikzlibrary{fadings}
\usetikzlibrary{patterns}
\usetikzlibrary{shadows.blur}
\usetikzlibrary{shapes}
\usepackage{float}

\newcommand{\nn}{\nonumber}
\newcommand{\del}{\delta}

\title{Optimal Lockdown Strategy in a Pandemic: An Exploratory Analysis for Covid-19\footnote{This research did not receive any special grant from any funding agencies in the public, commercial, or not for profit sectors. Further, none of the authors has any relevant material or financial interests that relate to the research described in this paper.}}

\author{Gopal K. Basak\thanks{
Stat-Math Unit, Indian Statistical Institute, Kolkata 700108, India, Email: gkb@isical.ac.in . \hspace{5pt} The author wishes to thank Vivek Borkar for his initial communication and fruitful discussion between them on SIAR model that includes asymptomatic individuals, raising the author's research interest in this direction.}, Chandramauli Chakraborty\thanks{Indian Statistical Institute, Kolkata 700108, India, Email: chandramaulichak108@gmail.com}, 
Pranab Kumar Das\thanks{
Centre for Studies in Social Sciences, Calcutta, R1, B.P. Township, Kolkata 700094, India, pkdas@cssscal.org
}}

\date{}
\begin{document}

\maketitle

\vspace{-32pt}
\abstract{The paper addresses the question of lives versus livelihood in an SIRD model augmented with a macroeconomic structure. The constraints on the availability of health facilities - both infrastructure and health workers determine the probability of receiving treatment which is found to be higher for the patients with severe infection than the patients with mild infection for the specific parametric configuration of the paper. Distinguishing between two types of direct intervention policy - hard lockdown and soft lockdown, the study derives alternative policy options available to the government. The study further indicates that the soft lockdown policy is optimal from a public policy perspective under the specific parametric configuration considered in this paper.

}

\bigskip


\medskip

\section{Introduction}
The objective of the present paper has been evaluation of the policy of direct intervention to arrest the spread of the Covid-19 pandemic. Using an aggregative framework the paper incorporates the evolution of the infected population from asymptomatic to mild and then to severe depending on the nature of required treatment and highlights the role of binding health infrastructure. Model simulations for plausible parameter values show that a policy of no intervention finds the pandemic to end after around 630 days, when herd immunity is reached, but with a substantial loss of lives as well as more than 1\% fall in the aggregate output. If, however, a strict lockdown policy is implemented as and when the health constraint binds, the loss of lives is only about 0.77\% for the same period compared to no lockdown case, but the economic loss is about 5\%. A soft lockdown policy that lifts the extent of the lockdown depending on how far the health infrastructure is away from the binding level produces outcomes in between these two extremes in terms of both loss of lives and economic loss. This model should not considered as a forecasting model, rather it aims to explore the evolution of the disease over time. The parameter values have been chosen keeping a country like India in mind. However, these can be amended to characterize any other country with different institutional structure.

The epidemiology model in economics is not entirely new, but there has been a spurt in the volume of the published papers and work in progress with the onset of Covid-19.\footnote{The classic reference on epidemiology model is Kermack and McKendrickís (1927). More recent literature includes Hethcote (2000), Chowell {\em et al.}(2009) among others. Anastassopoulou {\em et al.} (2020), Bertorzzi {\em et al.} (2020) and Sameni (2020) are a few of the epidemiology models addressing specifically Covid-19 pandemic.}  This literature employs what can be called macro SIR model (or some of its variant) that supplements the epidemiology laws governing the spread of the virus with a set of additional equations to reflect the economic behavior in a representative agent framework.\footnote{Brodeur {\em et al.} (2021) provides a good survey on the literature of economic issues of Covid-19 covering various aspects of which macro SIR models form a part.} The spread of virus gives a negative shock to the supply of labour as infection starts rising leading to death of a large population. This in turn transmits negative shocks to aggregate production, consumption etc. and the subsequent rounds leading to what has often been called pandemic led recession (Gregory {\em et al.}, 2020; Guerreri {\em et al.}, 2020; Kaplan {\em et al.} 2020). The process continues until a cure, usually a vaccine becomes available or when the herd immunity is reached transforming the epidemic or the pandemic into an endemic phenomenon. 

The justification for intervention policy stems from the fact that there is an externality in an epidemic (and more so in a pandemic) as uninfected individuals become susceptible to the disease when they come in contact with infected individuals in the work place or any other gathering (Farboodi {\em et al.}, 2020; Bryant {\em et al.}, 2020; Dimdore-Miles {\em et al.}, 2020). Thus the need for public action has been advocated in the literature in various forms. Bethune and Korinek (2020), Berger {\em et al.} (2020), Gatto {\em et al.} (2020), Grigorieva {\em et al.} (2020) suggest containment of infected people so as to reduce the rate of infection. Eichenbaum {\em et al.} (2020) proposes a model with containment measure in the form of consumption tax and a lump sum transfer. A tax on consumption tax reduces it and also makes leisure more attractive so that mixing of people both for purchase of goods and supply of labour decrease so that infection rate gets reduced. The tax is rebated to households so that disposable income remains unchanged. Jones, Philippon, and Venkateswaran (2020) argues the case for social distancing. Direct lockdown that suspends economic activities and other gatherings is favoured by Alvarez {\em et al.} (2020), Aspri {\em et al.} (2021), Caulkins {\em et al.} (2021), Palma, {\em et al.} (2020), Rawson {\em et al.} (2020). 

Lockdown policy to arrest spread of infection has a wedge between loss of lives versus loss of livelihood. Hence, Kaplan {\em et al.} (2020) suggests fiscal stimulus during the lockdown period to compensate for economic loss. Rampini {\em et al.} (2020) has argued the case for sequential lockdown. However, Born {\em et al.} (2020) suspects efficacy of lockdown policy. Using age specific demographic profile of population Acemoglu {\em et al.} (2020) shows that a targeted lockdown has a lower cost in terms of lower GDP as well as lower fatality rate, while Gollier {\em et al.} (2020) questions the efficacy of premature withdrawal of lockdown. With this introduction we organize rest of the paper into 3 sections – Section 2 proposes the model, Section 3 the simulation results and discussion and the last section concludes with scope of future work.

\section{Model}
The model in this paper follows a standard SIRD model for the evolution of the disease dynamics, but employs a more disaggregated framework with two kinds of infections, viz. mild infection and severe infection. In the former case the patients do not need hospitalisation while in the second case hospitalisation is very much required. Hence, two crucial factors that are very important in our analysis for reducing the death rate are the availability of doctors including other health workers in the former case and both doctors (including health workers) and hospital beds for treating the severely infected patients. It has been evidenced across the globe that health infrastructure, viz. both doctors and hospital beds become very important determinants in curbing the spread of infection and hence death. With the rise in the number of either or both of the mild  or  severely infected patients to a very high level available doctors cannot provide treatment and/ or hospital beds become unavailable. In this situation imposition of lockdown becomes imminent to arrest the number of infected patients. Imposition of lockdown can also be justified to reduce the mixing of population.     

The total population is assumed to remain fixed at $N$ during the period of analysis, but there are two kinds of susceptible population – general, $S_{g,t}$ and health workers, $S_{h,t}$ at any $t$; the latter includes both doctors and other health workers. It is assumed that there are three stages of infection - asymptomatic, mild and severe. In the beginning the infected population is assumed to remain asymptomatic ($A_t$) for a few days, then some of them recover, the rest , $I_{m,t}$ starts showing mild symptoms. They need treatment in the form of OPD or telephonic advice, but do not require hospitalisation. Some of themreceive treatment and the rest do not  depending on whether enough doctors are available or not. In either case a proportion of the mildly infected population becomes severely ill ($I_{c,t}$) and the rest recover with a higher probability if receive treatment than if does not. The severely ill pool of patients receive treatment and the rest do not receive treatment. In this case availability of treatment can be constrained by either of the availability of doctors or the availability of hospital beds.
 
The evolution of the disease dynamics described above are presented in terms of the following equations.
\begin{align}
\label{eq1}
 & \Delta{S_{g,t+1}} = -\lambda_{g,t} S_{g,t} A_t/N\\
 & 
 \Delta{S_{h,t+1}} = -\lambda_{h,t} S_{h,t} (\lambda_{0,t} A_t + I_{c,t} + \lambda_{1,t} I_{m,t})/N
 \\
 & 
 \Delta{A_{t+1}} = \lambda_{g,t} S_{g,t} A_t/N + \lambda_{h,t} S_{h,t} (\lambda_{0,t} A_t + \lambda_{1,t} I_{m,t} +I_{c,t})/N 
\nn\\
 & \hspace{2.0cm}
 - \beta_{0} \alpha_{0} A_t - \beta_{0} (1-\alpha_{0}) A_t
 \\
  &
  \Delta{I_{m,t+1}} = \beta_{0} \alpha_{0} A_t - \beta_1(\alpha_{m,t}(1-\alpha_{22})+(1-\alpha_{m,t})(1-\alpha_1))I_{m,t}
  \nn\\
  & \hspace{2.0cm}
  - \beta_1(1-\alpha_{m,t}(1-\alpha_{22})-(1-\alpha_{m,t})(1-\alpha_1))I_{m,t}
  \\
 & 
 \Delta{I_{c,t+1}} = \beta_1(1-\alpha_{m,t}(1-\alpha_{22})-(1-\alpha_{m,t})(1-\alpha_1))I_{m,t}
 \nn\\
  & \hspace{2.2cm}
  -\gamma_1(\alpha_{c,t}(1-\alpha_{42})+(1-\alpha_{c,t})(1-\alpha_3))I_{c,t}
  \nn\\
  & \hspace{2.4cm}
  -\gamma_1(1-\alpha_{c,t}(1-\alpha_{42})-(1-\alpha_{c,t})(1-\alpha_3))I_{c,t} 
  \\
 & 
 \Delta{R_{t+1}} = \beta_{0} (1-\alpha_{0}) A_t + \beta_1(\alpha_{m,t}(1-\alpha_{22})+(1-\alpha_{m,t})(1-\alpha_1))I_{m,t}
 \nn\\
 & \hspace{2.0cm}
 + \gamma_1(\alpha_{c,t}(1-\alpha_{42})+(1-\alpha_{c,t})(1-\alpha_3))I_{c,t} 
 \\ 
  & 
  \Delta{D_{t+1}} = \gamma_1(1-\alpha_{c,t}(1-\alpha_{42})-(1-\alpha_{c,t})(1-\alpha_3))I_{c,t}
  \\
\label{eq12}
 & (S_{g,t}+S_{h,t})+A_t+(I_{m,t}+I_{c,t})+R_t+ D_t = N
\end{align}
where, $\lambda_{g,t}$ = rate of change of susceptible general population to asymptomatic population over time when mixing with the  asymptomatic population, \\
$\lambda_{h,t}$ = rate of change of susceptible health workers to asymptomatic population over time when mixing with the critical patients being treated and proportion of mild patients being treated and a proportion of asymptomatic population, \\
$\lambda_{0,t}$ = proportion of asymptomatic population who come in contact with the health care workers, \\
$\lambda_{1,t}$ = proportion of the mildly infected population under treatment who come in contact with the health care workers, where
$\lambda_{0,t}$ and $\lambda_{1,t}$ change their values only when lockdown is imposed or withdrawn,
\\
$\beta_0$ = rate of depletion from the pool of asymptomatic population, \\
$\alpha_0$ = probability that an asymptomatic person gets mildly infected (as opposed to getting recovered), \\
$\alpha_{m,t}$ = probability that a mildly infected patient can not receive treatment at time $t$, \\
$\beta_1$ = rate of depletion from the mildly infected population under treatment, \\
$\alpha_{1}$ = probability that mildly infected  patient under treatment falling critically ill, \\
$\alpha_{22}$ = probability that a mildly infected untreated patient falling critically ill, \\
$\alpha_{c,t}$ = probability that a critically ill patient does not get treatment at time $t$, \\
$\alpha_{42}$ = probability of death of untreated critically ill patients, \\
$\gamma_1$ = rate of depletion from the treated critically ill population, \\
$\alpha_{3}$ = probability of death of treated critically ill population. \\
The disease dynamics is represented in terms of a flow chart in 
figure 1.

\bigskip
Fig. 1 is about here.
\bigskip

It may be noted that $\lambda_{g,t}$ is the outcome of two parameters, viz. rate of mixing of the population, $c_m$, and rate of infection of the disease, $c_i$. These two can be combined to generate $\lambda_{g,t}$ as in below:
\begin{equation}
\label{eq13}
\lambda_{g,t} = c_m c_i \theta_{t}^{1+\nu}
\end{equation}
where $\nu$ is a parameter and $\theta_{t}$ is lockdown parameter, such that $0<\theta_t <1$. Similarly for the $\lambda_h$ term. With no lockdown $\theta_t$ takes the value of unity so that the rate of transmission, $\lambda_{g,t}$ has full impact. With $\theta_t<1$, which happens when a policy of lockdown is implemented, the rate of transmission $\lambda_{g,t}$ gets reduced.  
The rates of change, viz. $\lambda_{g,t}$, $\lambda_{h,t}$, $\lambda_{0,t}$, $\lambda_{1,t}$, and the probability terms (transition probabilities), viz. $\alpha_{m,t}$ and $\alpha_{c,t}$, $\alpha_{41,t}$ are indexed by time implying that these vary with time. The rates of change terms vary depending upon whether there is lockdown or not and if there is lockdown then type of lockdown, viz. hard or soft. The probability of getting treatment vary over time because of the availability of doctors and beds which in turn are dependent on policy opted for lockdown. 

Denoting the number of available doctors for each of the mild and severely infected patients (doctor-patient ratio) by $\phi_{m,t}$ and $\phi_{c,t}$ respectively on each day and total beds available for treatment by $B$ the time dependent rates of change and the transition probabilities are defined below: \\
\begin{equation}
\label{eq16}
\alpha_{m,t} 
= 1-\frac{exp(\lambda_m/(1 - min(h_{g,t}, 1)))}{exp(- \lambda_m)}
\end{equation}
and
\begin{eqnarray}
\label{eq17}
&&
\alpha_{c,t} 
= 1-\frac{exp( - max(\lambda_c/(1 - min(h_{g,t}, 1)),  \lambda_b/(1 - min(h_{b,t}, 1))))}{exp(- max(\lambda_c, \lambda_b))}  \\
&& \mbox{where } \ h_{b,t} = \frac{I_{c,t}}{B_t}, \ 
h_{g,t} = \frac{(I_{m,t} \phi_{m,t} + I_{c,t} \phi_{c,t})}{(S_{h,t} + R_{h,t-\delta_h}+A_{h,t})} 
\nn
\end{eqnarray}
with $\phi_{m,t} = 1/15$ (if there is no lockdown),   $\phi_{m,t} = 1/26$  (if there is lockdown) and $\phi_{c,t} = 1/7$  (if there is no lockdown),   $\phi_{c,t} = 1/10$  (if there is lockdown), \\
$(S_{h,t} + R_{h,t-\delta_h}
+A_{h,t}
) = $ no. of available doctors at time $t$.
The no. of beds, $B$ at time $t$ (assumed constant, but in reality, however, during the stressed time other facilities are converted into beds/infrastructure),

\bigskip

For the sake of simplicity the services of the doctors and other health workers have been clubbed together in the same category. Separate constraints on the availability of each of the other health workers can be introduced in a more generalized and real life model. However, the constraint on the availability of doctors can capture other health workers also when there is a fixed proportion of requirement between doctors and other health workers. The availability of beds can also be thought of as a composite good that includes other peripherals like facilities of ICU, availability of oxygen and other medicines etc. In a situation when the disease spread reaches very high level, either of the availability of doctors or that of the beds binds. However, the available data shows that it is the beds that becomes most critical for the treatment of the patients in the Indian context.

The labour force at any time point $t$ is the total population less the infected and dead plus recovered with the necessary delay to return to work, which may be said to be the adjusted sum of $S_{g}$ and $S_{h}$ (see \eqref{eq18}). However, the recovered has a delay of a few days to join labour force. In this respect we are making a distinction in respect of labour force for hospital services and other economic activities. The former is not subject to any lockdown. So lockdown is imposed on the labour force for general economic activities. Accordingly the general labour force is derived from \eqref{eq1} through \eqref{eq12} to be:
\begin{equation}
    \label{eq18}
L_t (\theta_t) = \theta_t(S_{g,t} + A_{g,t} + R_{g,(t-\del_{g})}) + a(S_{h,t} + A_{h,t} + R_{h,(t-\del_{h})})
\end{equation}
where $a$ is a multiplier used for the conversion of the output of a health worker to that of a general worker, $\del_h$ and $\del_g$ are the delay time for the corresponding worker to get back to the work and $\theta_t$ is the lockdown factor, $\epsilon (0,1)$. A value of $\theta_t$ = 1 implies no lockdown in force and all of the available labour force is allowed to work while imposition of lockdown means $\theta_t < 1$ so that only a part of the labour force is allowed to work. 

The policy instrument for imposition of lockdown in this paper is linked with the availability of doctors and hospital beds (vis-\'{a}-vis infrastructure). As and when either or both the constraints bind or close to the binding level lockdown is imposed via restriction on the use of labour. It remains in place so long as the constraint is relaxed to a reasonably low level. This is called {\it hard lockdown} in this paper. However, in order to functioning of the essential services, such as shopping for food articles, electricity, gas, banking, health etc. a certain minimum level of activities is allowed when lockdown is in force. The lockdown function takes two values and is defined at each $t$ as in below:
\begin{align}
\label{eq19}
& \theta_t = \theta_0 \ \mbox{ if } h \ge 1   \nn\\
& \hspace{0.3cm} = 1 \ \mbox{ if } l <  1
\end{align}
where $\theta_0$ is the level of activities allowed under lockdown, and $h$ and $l$ refers to the maximum and normal levels of of stress that health facilities can accommodated and defined as in below:
\begin{align}
 \label{eq20}
 h = Min(S_{h,t} + A_{h,t} + R_{h,t},B)
 \end{align}
\begin{align}
    \label{eq21}
l= Min(\frac{2}{3} B, (S_{h,t} + A_{h,t} + R_{h,t})/\sqrt{\frac{26\times10}{7\times15}}) 
\end{align}
There is another option for the lockdown policy, viz. the degree of lockdown is set to be determined by a rule depending on the extent of the constraints on the availability of doctors and / or hospital beds. So, $\theta_t$ is now replaced by a time dependent function. This policy of lockdown has been called {\it soft lockdown} in this paper. The rule for changing the extent of lockdown can mainly be thought of three kinds, viz. (a) a constant proportion of the deviation of actual availability from the binding level, called linear, (b) an increasing function of the deviation from binding level, called convex, and (c) a decreasing function of the deviation from binding level, called concave, defined as in below:
\begin{align}
\label{eq22}
\mbox{(a) } \ \ \theta_t 
& = \theta_0 \frac{\kappa - l}{h - l} + \frac{h - \kappa}{h - l}  \nn\\
& = (\theta_0 - 1) \frac{\kappa - l}{h - l} + 1   \nn\\
\mbox{(b) } \ \ \theta_t
& =  (\theta_0 - 1) f_{conv}\left(\frac{\kappa - l}{h - l}\right) + 1  \nn\\  
\mbox{(c) } \ \ \theta_t
& =  (\theta_0 - 1) f_{conc}\left(\frac{\kappa - l}{h - l}\right) + 1  \nn\\  
\end{align}
where $f_{conv}, \  f_{conc}$ are convex and concave increasing  bijections from $[0, 1] \mapsto [0, 1]$, respectively and 
\begin{align}
\label{eq22a}
\kappa = Min(\frac{I_{c,t}}{7}+\frac{I_{m,t}}{15},\frac{3}{2}I_{c,t})
\end{align}

A power function of the form $Z^\mu$ with $0<\mu < 1$ represnts a concave adjustment rule for $\theta$, with $\mu>$1 represents a convex adjustment rule and a $\mu=1$ represents a linear adjustment. A value of $\mu$ equal to zero represents hard lockdown with $\theta$ equaling $\theta_0$.

The economic activities, measured by GDP in this model, is represented by an aggregate Neo-classical production function of Cobb-Douglas variety with constant returns to scale in two factors employed at time $t$ - capital , $K_t$ and labour$L_t$ and a technology parameter $V$. With capital, $K_t$ and the technology parameter, $V$ remaining constant in short period of time the product $V K^{1-\alpha}_t$ is normalized at unity.Thus aggregate output, $Y_t$ is given by, 
\begin{equation}
\label{eq23}
Y_t= L^{\alpha}_t
\end{equation}
where $0< \alpha <1$. 
Any variation in aggregate output is obtained by varying labour. The benchmark full employment level output, $\Bar{Y}$ is the output corresponding to full employment of labour $\bar{L}$ before the pandemic began, i.e. $N$. It is expected that actual output during the pandemic period is well below $\bar{L}$ level, even when there is no lockdown, because population growth is ignored in the model while the labour force decreases due to death. In this model total population is considered to be the working population or the labour force in the age group 15 to 64. That is, the population, $N$, is equated with the entire labour force for the sake of simplicity. This age group is around 65\% of total population in India. 

We do not distinguish between different components of GDP for two reasons. First, unlike Bethune and Korinek (2020), {\em et al.} Eichenbaum {\em et al.} (2020), Hall {\em et al.} (2020), Kaplan {\em et al.} (2020), we do not aim to find the optimal lockdown policy as a central planner's welfare maximisation problem. Instead, the optimum lockdown in this model is obtained by minimizing a loss function defined in terms of deviation of aggregate output and number of deaths from pre-specified targets, viz. pre-lockdown levels.  
This approach does not need any distinction between components of aggregate output, all that matters is the aggregate output. Secondly, unlike Jones {\em et al.} (2020) and others, we do not distinguish between meeting of people on the basis of consumption and production purposes.

When lockdown (hard or soft lockdown) is in place the aggregate output is not given by \eqref{eq23}, instead it is amended to incorporate the lockdown policy parameter $\theta_t$ with a policy of lockdown  (either of the hard or soft) in place as in below:
\begin{equation}
\label{eq24}
Y_{l,t} = (L_t(\theta_t))^{\alpha}
\end{equation}

Finally, the loss function, $\Psi$ for the evaluation of the policy performance is defined in terms of the deviation of GDP, $Y_t$ from a target level, $\Bar{Y}$ and total fatality, $D_t$ from a target level, which in this case is assumed to be zero over the usual death from other diseases. As these two items are not conformable for addition we use a weight, $\chi$ which represents the statistical value of life. The objective function for the policy makers is given by \\ 
\begin{align}
\label{eq25}
& \Psi_t = \sum_{j=0}^{t_0}  [(1- \frac{Y_{l,t+j}(\theta_t)}{\bar{Y}})\bar{Y} )^m +\chi (D_{t+j}(\theta_t))^m ]
\end{align}
where, $t_0$ is the terminal period for which the policy of lockdown is implemented. The control of the minimization exercise is $\theta_t$, the type and extent of lockdown under alternative policy regimes and $Y_{l,t}$ and $D_t$ are both state variables. The parameter $m$ is the power of the loss function. For $m$ = 1, it is linear in target output and number of deaths, though convex in the control $\theta_t$. For $m$ =2, it becomes the standard quadratic loss function in target output and number of deaths. It may also be noted that $m$ is an argument in the $\Psi$ function. So minimisation of $\Psi$ by choice of the control $\theta$ implies the choice for a given $m$. 

This form of the objective function has important implications. When the (policy) measures taken fails to reach the targets then exponential (linear) loss penalizes severely compared to, say, quadratic loss, whereas when the targets are achieved (or exceeds) benefit is less compared to that of quadratic loss.

The spread and control of a pandemic in a short horizon problem, in the specific case  of Covid-19 is expected to continue for 3 years if Spanish flu is an indicator. Thus the objective function is not discounted for the control of a pandemic is a short horizon problem. Statistical value of life $\chi$ is drawn from the standard literature (Shelling, 1968) defined as an estimate of the financial value that society places on reducing the average number of deaths by one. Employing the widely adopted method of Viscusi and Kinpand (2003) the statistical value of life is estimated to be around USD 0.43 million per death at current prices in India. Majumder and Madheswaran (2018) provided a more liberal estimate to the tune of USD 0.64 million for the Indian population. This estimate is used in the present paper.

Thus the problem of the policy maker is given by:
\begin{equation}
    \label{eq26}
    \min_{\theta_t} \ \Psi_t
\end{equation}
s.t. \eqref{eq1} to \eqref{eq12}, 
\eqref{eq16}, \eqref{eq17}, 
\eqref{eq19} (or \eqref{eq22} depending upon lockdown regime type), \eqref{eq24}.  

This completes the description of the model economy with the disease dynamics. In the next section when we undertake the simulation of the model we adopt two approaches to explore the policy options for the control of the spread of infection. First, we will take up a general characterization of the intertemporal trajectory with or without lockdown policy. It lists different levels of aggregate output and number of deaths with or without lockdown policy. Second, we provide the optimal bundle of aggregate output and number of deaths obtained via optimization of the loss function.

\section{Simulation}
The simulation to analyze the evolution of the disease over time is undertaken with plausible parameter and initial values are provided in Table 1.
The rates and probabilities, viz. $\lambda$ ‘s, $\alpha$ ‘s and $\beta$ ‘s are calculated using available data from governmental and international sources (www.mohfw.gov.in , www.worldometers.info/coronavirus etc.) . These are further compared with the values used in other studies for India and other countries to arrive at reasonable and meaningful estimates. The parameter of the aggregate production function, $\alpha$ is the share of wages in GDP, which is obtained from ILO (2018).  As already explained in the previous section the model in this paper considers the working population at the age group 15 to 64, which in the Indian case constitutes around 65\% of total population of 1.37 Billion. The number of health workers (doctors in this paper) is obtained from the World Development Indicators (WDI) database of the World Bank duly adjusted for the 65\% of the population. The number of hospital beds is also obtained from the WDI database duly adjusted in the same way.  

As for the parameters $\phi_{m,t}$, $\phi_{c,t}$ in \eqref{eq16} and \eqref{eq17} we consider  reasonable values in the twin cases when there is no lockdown and when  there is lockdown. In the former case the health system is under normal level of operation with a lower probability of not getting treatment while in the latter case the health system is stressed because of a very high level of the spread of the disease with too many patients and the probability of not getting treatment rises. The form of the probability functions \eqref{eq16} and \eqref{eq17} ensures that with rising number of patients there is stress on the health system and the probabilities of no treatment increases. We further introduce a mutation of the virus after $450^{th}$ day, i.e. after five quarters from the time of occurrence of the disease. This is introduced by an exogenous increase of $40$\% in the rate of transmission of infection, $c_i$ for the general population as well as among the health workers. This leads to an increase of $\lambda_{g,t}$ and $\lambda_{h,t}$ from $0.12$ to $0.168$ for the general population and from $0.08$ to $0.1$ for the health workers respectively in the absence of lockdown (i.e. $\theta_t$=1) as indicated in the footnote to 
Table 1.

However, it may be noted that the value of $\lambda_{g,t}$ and $\lambda_{h,t}$ decreases  when lockdown (hard or soft) remains effective, because during lockdown there is restriction on the mixing of the population given by $\theta_t$ and hence the disease spreads at a lower rate. 

\bigskip

                  Table 1 is about here
\bigskip

First, we consider the case in the absence of any lockdown policy. In this case the virus spread is allowed to take place without any check. The total number of infections (of the twin categories), recovery and death and proportion of fatality for the general population as well as for the health workers (measured in terms of observed recoveries) are provided in Table 2 
at the end of each quarter. The evolution of the disease is also presented in figure 2 through figure 6. 
The disease spreads at a high rate until the third quarter, thereafter slows down, reaches its peak between $4^{th}$ and $5^{th}$ quarter, shows sign of decline, but again increases with the occurrence of the new variant with higher rate of infection ($\lambda_g$ and $\lambda_h$) after $5^{th}$ quarter ($450^{th}$ day). However, it reaches the peak soon, at around $500^{th}$ day, and then starts declining. It is evident from the figures as well as the tables that the spread of the disease stabilizes after $630^{th}$ day, i.e. after around quarter to two years. Thence it becomes endemic with total death a little more than $5.29$ million which is 0.598\% of the population. 

\bigskip

                Table 2 is about here
\bigskip

		 Fig 2 is about here
		 
		 Fig 3 is about here
\bigskip

The fatality rate is 2.81\% as on the $630^{th}$ day. With the infection at the level of 53.46\% of the population the herd immunity is reached after $630^{th}$ day. \footnote{The herd immunity is calculated by total asymptomatic cases as percentage of starting population. } It may be noted that the herd immunity is reached with the spread of the disease among $53.6$\% of the population which is well below the expected value of above $60$\% as the forecasting models of the epidemiology literature predicts. It so happens in this model because we have considered only one type of population with the same rate of mixing. The rate of mixing differs among working population and non-working population and the children. Hence a more realistic model with realistic demography will take care of this apparent anomaly.

\bigskip

		Fig. 4 is about here
		
		Fig. 5 is about here
		
		Fig. 6 is about here
\bigskip

The probability of no treatment for both the mild as well as severe infections rapidly rises with the spread of infection from $3^{rd}$ quarter which attains very close to unity immediately. The probabilities in both the cases remain at unity until the infection starts declining after around $500^{th}$ day. However, the probability of no treatment for the mild case shows a few small dips.  Once herd immunity is reached there is no significant change in the total number of asymptomatic, mild or severe cases of infection or the number of death.  If another mutation with higher infection rate occurs then the date of stabilization would further extend but may not be too far, because with rising infection of the population number of cases will increase sufficiently to reach the level of herd immunity. The GDP falls about 1\% in this case reaching the lowest point on $630^{th}$ day, there after starts recovering. The pace of recovery is slower than the fall reflected in a lower (absolute value) slope of the former than the latter. 

This case is hypothetical as the no lockdown policy was not in force in India or elsewhere across the globe \footnote{Sweden was the only country which did not impose lockdown initially, but they implemented other measures to arrest spread pf the disease.}, but it is a benchmark for comparing different policy options. However, the relevance of this case stems from the fact that it provides a comparison with the Spanish flu of 1918-20. The Spanish flu continued for about three years but the cases of death occurred mainly in the first year across the globe with 1.42\% of population in the first year and 2.1\% for the three years period taken together. The corresponding rates were 4.1\% and 5.22\% respectively for India which was the second highest (after Kenya). The three year aggregate rate for European countries and USA were well below less than 1\%. The medical infrastructure at that time was much less by today’s standard and for India it was almost non-existent for the native population. Though quarantine was implemented in India and elsewhere, no large scale lockdown was imposed. So largely it is the herd immunity that arrested the spread of the disease. 

Next we consider the simulation with a policy of hard lockdown presented in Table 3
and figures 7 through 11.
As was discussed earlier hard lockdown is implemented by the rule given by \eqref{eq19}, i.e. as and when either of the demand for doctors or hospital beds exceeds the availability and it is lifted once the constraint is relaxed. Accordingly the hard lockdown policy is implemented (lifted) on the $306^{th}$ ($336^{th}$), $463^{th}$ ($504^{th}$), $551^{th}$ ($593^{th}$), $645^{th}$ ($684^{th}$), $744^{th}$ ($782^{th}$) day. We have shown the case of hard lockdown policy with $\theta_t$ =0.5, i.e. allowing $50$\% workforce. Of the total workforce $30$\% comprises of the employment in the essential services, such as health care, transport etc., the rest from the remaining $70$\%. It is evident from figure 7
that as the disease starts rising very fast after $3$rd quarter as in the no lockdown case leading to constraint on the health services. Then the policy of lockdown is imposed from $306^{th}$ day. This 

\bigskip				
                              Table 3 is about here
\bigskip

				Fig. 7 is about here
				
				Fig. 8 is about here
\bigskip

\noindent
reduces the mixing of the population and hence the value of $\lambda_g$ (and also $\lambda_h$) leading to relaxing the constraint as given by \eqref{eq19}. The lockdown is then withdrawn after about a month on the $336^{th}$ day. Withdrawal of lockdown allows the spread of the disease at the previous rate of $\lambda_g$ until it reaches a level for \eqref{eq19} to bind and the process of imposition of lockdown repeats itself. With the arrest of the spread of the disease temporarily during the period of lockdown the hospital services gets better off but it as well reduces the process of herd immunity and hence the steady state. In fact only $13.7$\% of the population attains herd immunity on the $630$th day and $20$\% on $810$th  day in this case. The herd immunity is reached well after $2200$ days (not shown here in figure). The total number of infections as well as fatalities is much lower with hard lockdown until $7^{th}$ or $9^{th}$ quarter. The total cases of fatality on the $630^{th}$ and $810^{th}$ day are $0.77$\% and $0.122$\% of the population respectively. 

\bigskip

				Fig 9 is about here
				
				Fig 10 is about here
				
 				Fig 11 is about here
\bigskip

The probability of no treatment for both types of patients show different pattern in this case compared to the case of no lockdown. It is evident from figures 9 and 10
that probability of no treatment for the mild case never reaches unity but the same for the severe case remains close to unity most of the time. This happens because availability of hospital beds has stronger role as a constraining factor than the availability of doctors. Once the lockdown is imposed, $\alpha_{m,t}$ and $\alpha_{c,t}$ start falling. The policy of lockdown helps reduce number of infections as well as death but it continues in the form of a pandemic for longer time. However, the contraction of the GDP is around 5\% and 6\% until eight and nine quarters respectively which are higher than in the no lockdown case. A comparison with the policy of no lockdown case highlights the issue of lives versus livelihood. With no lockdown death is much higher after about one and half years when the disease reaches steady state but the contraction in the GDP is lower and the economy begins early recovery. On the other hand a policy of hard lockdown reduces infection and death in current times but extends the period of pandemic and loss of output is much higher.

Now we consider the case of soft lockdown. In the previous section we proposed three cases of soft lockdown based on policy rule pertaining to the constraints on the availability of doctors and hospital infrastructure, viz. linear, convex and concave, assuming that the implementation started at some level of stress in the health system in comparison to that of the hard lockdown. Tables 4 and 5 
provide the simulation results for convex and concave rules respectively, assuming it was initiated at 75\% stress level of the health system. Both the cases show that the soft lockdown policy lies intermediate between the policy of no lockdown and the policy of hard lockdown in terms of levels of infection and case fatality. However, infection and fatality are 

\bigskip

				Table 4 is about here
\bigskip

				Fig. 12 is about here
				
				Fig. 13 is about here
				
				Fig. 14 is about here
				
				Fig. 15 is about here
				
				Fig. 16 is about here
\bigskip

				Table 5 ia about here
\bigskip

				Fig. 17 is about here
				
				Fig. 18 is about here
				
				Fig. 19 is about here
				
				Fig. 20 is about here
				
				Fig. 21 is about here
\bigskip

\noindent
higher in the convex case than in the concave case. The fall in GDP is lower in the convex case than in the concave case. This happens because of the fact that in the concave rule the extent of the lockdown is stringent than in the convex rule. In other words the $\theta_t$ value is lower in the former than in the latter. The probability values of no treatment for both mild and severe infections follow similar pattern though somewhat higher in the convex rule. In fact the probability of no treatment for the severe case show persistence at unit value while in the concave rule it changes frequently.
The linear policy rule is intermediate between these two cases. Hence, we do not provide the result here. 

The above policy rules give options to the policy makers between lives versus livelihood. The policy of no lockdown has lower adverse impact on income and employment but a higher fatality in the current times while the opposite happens for the policy of hard lockdown. The policy of soft lockdown, which has been practiced in many contexts, looks attractive to the policy makers of a popular government, because it does not contract income and employment as much as for the hard lockdown rule but has lower fatality in the current times than in the no lockdown case. This choice problem assumes its importance in the absence of a cure in the form of vaccine or even if a vaccine is available, its efficacy in the face of mutation of the virus. It is upto the policy makers to decide about the lives versus livelihood question and accordingly takes decision about the implementation of the lockdown policy.

Finally, we provide an evaluation of the policy regimes in terms of minimization of the loss function $\Psi_t$, vide \eqref{eq25} over the period of analysis, viz. $810$ days. Among the various forms of the loss function based on the power $m$, we provide here only the case of loss function linear in the deviation of output from target and death. However, it is clear that it is convex in both the degree of lockdown, $\theta_t$ and adjustment rule for relaxing or tightening the it, $\mu$. 
Table 6 and figure 22
describe this case. As $\mu$ increases the policy regime moves more towards lesser degree of lockdown and vice versa. It is evident that the minimum cost is achieved for some $\theta$ in the interval (0.4  0.5) for $\mu$ = $1/100$. The latter implies that as lockdown tightens the $\Psi_t$ function approaches the minimum.   

\bigskip

			Table 6 is about here
\bigskip

			Fig. 22 is about here
\bigskip

As $\mu$ is raised the effect of loss of GDP more than outweighs the cost of death, because in the loss function the cost of life does not increase proportionately. Hence as $\mu$ rises, say for $\mu$ =$2$ implying a quadratic loss function in the deviation of target output and death gives choice for no lockdown. On the other hand a higher value of life, $\chi$ may give higher importance to life, but unlikely to favour a policy for no lockdown or soft lockdown.

\section{Conclusion}

Pandemic like Covid-19 threw everyone, from scientists to politicians, common people, government worldwide in the same page. To find out ways to stop/reduce rate of infection and to reduce fatalities are need of the hour. Since the disease had no known cure by medication, only way to let people treat symptomatically to relieve and let heal (with or without other related severe illnesses/side effects) or die. Barring certain age group (especially children and young adults) since the disease did not discriminate between poor and rich, or nationality or special sects, health system of all the countries in the world stretched beyond their limit to address issues during the peak of the infection. Only measures that could be adopted to ease the health system and not fall apart is to reduce the mixing of the people. Hence the policy of lockdown and other restrictions in various ways to reduce the contacts across population. But all these measures adversely affected economic activities worldwide, perhaps much more than any other crisis before. 

Our aim in this paper was to explore, whether there is a way to minimise the effect of the lockdown on the economy, whether there is any optimal way to set the restrictions to reduce mixing, depending on the situation of each country, contingent on the specific economic structure and the relative value of human life. We divide these restrictions in 3 general ways: (i) no lockdown, (ii) hard lockdown and (iii) soft lockdown. 

Our findings show that in the countries where relative values of human life is higher it could go with harder form of the lockdown to reduce the stress in the health system and attain the optimal level in balancing the economic activities, whereas the countries where the value of human life is relatively lower would be better off with softer version of the lockdown, that is allowing certain level of economic activities varying according to the stress on the health system and/or infection / fatality rate, etc. Thus the economy would be functioning relatively at a higher level of activities than with a policy of hard lockdown. As a result people with low level of earnings, generally daily wage earners, would not be adversely affected, especially when they are not compensated by the government when unemployed. This is especially true in low income countries with no effective way of unemployment registration with very large informal/ unorganized sector.

We have restricted our attention to one mutation of the virus (or say, 2nd wave, due to mutation), one can do similar exercise for third or more waves. Although, mutation had been introduced here exogenously keeping the fatality rate same, it can be introduced as an endogenous random phenomenon and with lower fatality rate, as it has happened in several countries. Similarly, repeated infection has not been introduced in this paper (even after mutation), but it can easily be incorporated and enhance the chances of the 2nd and more waves and more slow rate of going to stability. 
Under the model conditions, stability is reached faster in the no lockdown case (630 days) but with more deaths in the current times.  Hard lockdown achieves stability at a later date (more than 2200 days) with lower fatality, but reduced level of output as well. The policy of soft lockdown with a loss function linear in the deviation of output from target and number of deaths can help countries opt for suitable way to run the economic activities at the same time relieving the stress on the health system. The resultant outcomes of fatality, aggregate output and steady state of the disease spread lie somewhere in between the policies of hard lockdown and no lockdown. 

Purpose of this model is not to be considered as a forecasting model, rather it is an exploratory model for understanding the spread and arrest of the pandemic like Covid-19 over time. 
The model is simulated with the parameter values that may resemble  a country like India. With suitable amendments it can be used to study the characteristics of the disease elsewhere.
Further, this paper has explored only homogenous single (one large common and one small health worker) working population. One can incorporate demographic composition with an elderly age group and children group in this model to make it a model for complete population. However, qualitative results will still hold, with children group having lower infection and lower fatality rates and elderly having higher infection and fatality rates and proportion of them in the population varying from country to country.

The model structure can easily be generalized in bringing interacton with other (one or more) population groups across regions or countries through trade / travel (i.e., mixing between districts, between states / provinces as happened in Europe, India and the USA and elsewhere across the globe). Such an interactive model will allow the policymakers to the study the local characteristic for any regions keeping interaction with neighbouring regions and help implement more effective measures/decisions. Further, one can incorporate vaccination with varying success rate, in a dynamic way as it is being done in various countries to allow it to study eventual stability under vaccination.

Overall this is only an exploratory model, which would allow several other component to make it a more realistic framework for the policymakers to device a decision-making mechanisms.

\section*{Figures}

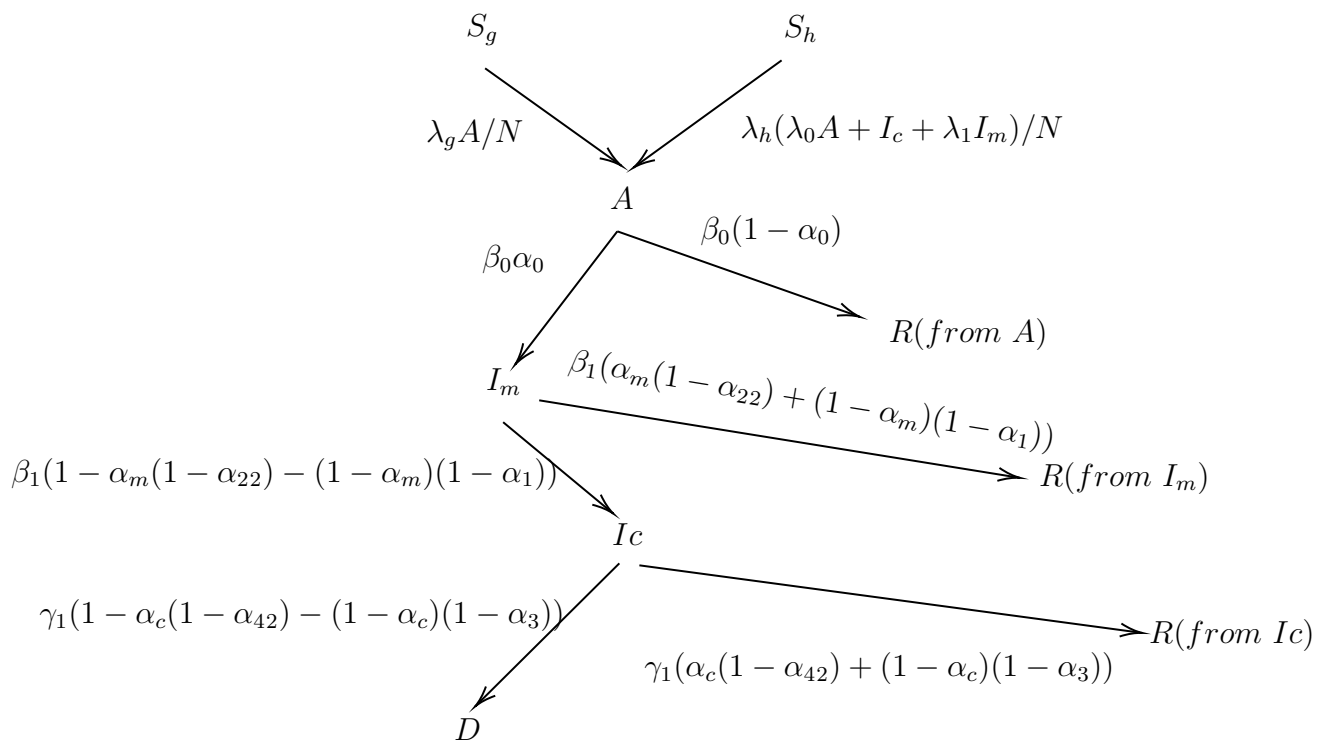
\begin{figure}[H]
\centering
\caption{Flowchart from susceptible to recovery / death}
\vspace{1cm}

\tikzset{every picture/.style={line width=0.75pt}} 

\begin{tikzpicture}[x=0.75pt,y=0.75pt,yscale=-1,xscale=1]

\draw    (246,59) -- (312.38,106.83) ;
\draw [shift={(314,108)}, rotate = 215.78] [color={rgb, 255:red, 0; green, 0; blue, 0 }  ][line width=0.75]    (10.93,-3.29) .. controls (6.95,-1.4) and (3.31,-0.3) .. (0,0) .. controls (3.31,0.3) and (6.95,1.4) .. (10.93,3.29)   ;
\draw    (394,55) -- (321.62,107.82) ;
\draw [shift={(320,109)}, rotate = 323.88] [color={rgb, 255:red, 0; green, 0; blue, 0 }  ][line width=0.75]    (10.93,-3.29) .. controls (6.95,-1.4) and (3.31,-0.3) .. (0,0) .. controls (3.31,0.3) and (6.95,1.4) .. (10.93,3.29)   ;
\draw    (312,141) -- (431.12,183.33) ;
\draw [shift={(433,184)}, rotate = 199.56] [color={rgb, 255:red, 0; green, 0; blue, 0 }  ][line width=0.75]    (10.93,-3.29) .. controls (6.95,-1.4) and (3.31,-0.3) .. (0,0) .. controls (3.31,0.3) and (6.95,1.4) .. (10.93,3.29)   ;
\draw    (312,141) -- (262.21,206.41) ;
\draw [shift={(261,208)}, rotate = 307.28] [color={rgb, 255:red, 0; green, 0; blue, 0 }  ][line width=0.75]    (10.93,-3.29) .. controls (6.95,-1.4) and (3.31,-0.3) .. (0,0) .. controls (3.31,0.3) and (6.95,1.4) .. (10.93,3.29)   ;
\draw    (255,237) -- (308.47,281.72) ;
\draw [shift={(310,283)}, rotate = 219.91] [color={rgb, 255:red, 0; green, 0; blue, 0 }  ][line width=0.75]    (10.93,-3.29) .. controls (6.95,-1.4) and (3.31,-0.3) .. (0,0) .. controls (3.31,0.3) and (6.95,1.4) .. (10.93,3.29)   ;
\draw    (323,309) -- (572.02,342.73) ;
\draw [shift={(574,343)}, rotate = 187.71] [color={rgb, 255:red, 0; green, 0; blue, 0 }  ][line width=0.75]    (10.93,-3.29) .. controls (6.95,-1.4) and (3.31,-0.3) .. (0,0) .. controls (3.31,0.3) and (6.95,1.4) .. (10.93,3.29)   ;
\draw    (313,308) -- (242.41,378.59) ;
\draw [shift={(241,380)}, rotate = 315] [color={rgb, 255:red, 0; green, 0; blue, 0 }  ][line width=0.75]    (10.93,-3.29) .. controls (6.95,-1.4) and (3.31,-0.3) .. (0,0) .. controls (3.31,0.3) and (6.95,1.4) .. (10.93,3.29)   ;
\draw    (273,226) -- (511.03,264.68) ;
\draw [shift={(513,265)}, rotate = 189.23] [color={rgb, 255:red, 0; green, 0; blue, 0 }  ][line width=0.75]    (10.93,-3.29) .. controls (6.95,-1.4) and (3.31,-0.3) .. (0,0) .. controls (3.31,0.3) and (6.95,1.4) .. (10.93,3.29)   ;

\draw (235,30.4) node [anchor=north west][inner sep=0.75pt]    {$S_{g}$};
\draw (387,28.4) node [anchor=north west][inner sep=0.75pt]    {$ \begin{array}{l}
S_{h}\\
\end{array}$};
\draw (307,117.4) node [anchor=north west][inner sep=0.75pt]    {$A$};
\draw (576,334.4) node [anchor=north west][inner sep=0.75pt]    {$R( from\ Ic)$};
\draw (245,208.4) node [anchor=north west][inner sep=0.75pt]    {$I_{m}$};
\draw (308,286.4) node [anchor=north west][inner sep=0.75pt]    {$Ic$};
\draw (229,384.4) node [anchor=north west][inner sep=0.75pt]    {$D$};
\draw (215,82.4) node [anchor=north west][inner sep=0.75pt]    {$\lambda _{g} A/N$};
\draw (372,80.4) node [anchor=north west][inner sep=0.75pt]    {$\lambda _{h}( \lambda _{0} A+I_{c} +\lambda _{1} I_{m}) /N\ $};
\draw (352,131.4) node [anchor=north west][inner sep=0.75pt]    {$\beta _{0}( 1-\alpha _{0})$};
\draw (243,146.4) node [anchor=north west][inner sep=0.75pt]    {$\beta _{0} \alpha _{0}$};
\draw (287.97,196.94) node [anchor=north west][inner sep=0.75pt]  [rotate=-9.27]  {$\beta _{1}( \alpha _{m}( 1-\alpha _{2}{}_{2}) +( 1-\alpha _{m})( 1-\alpha _{1}))$};
\draw (8,253.4) node [anchor=north west][inner sep=0.75pt]    {$\beta _{1}( 1-\alpha _{m}( 1-\alpha _{2}{}_{2}) -( 1-\alpha _{m})( 1-\alpha _{1}))$};
\draw (22,324.4) node [anchor=north west][inner sep=0.75pt]    {$\gamma _{1}( 1-\alpha _{c}( 1-\alpha _{4}{}_{2}) -( 1-\alpha _{c})( 1-\alpha _{3}))$};
\draw (324,351.4) node [anchor=north west][inner sep=0.75pt]    {$\gamma _{1}( \alpha _{c}( 1-\alpha _{4}{}_{2}) +( 1-\alpha _{c})( 1-\alpha _{3}))$};
\draw (446,183.4) node [anchor=north west][inner sep=0.75pt]    {$R( from\ A)$};
\draw (521,255.4) node [anchor=north west][inner sep=0.75pt]    {$R( from\ I_{m})$};

\end{tikzpicture}
\label{fig1}
\end{figure}

\begin{figure}[H]
\vskip-2cm
\hskip-2cm
\includegraphics[width=1.4\textwidth]{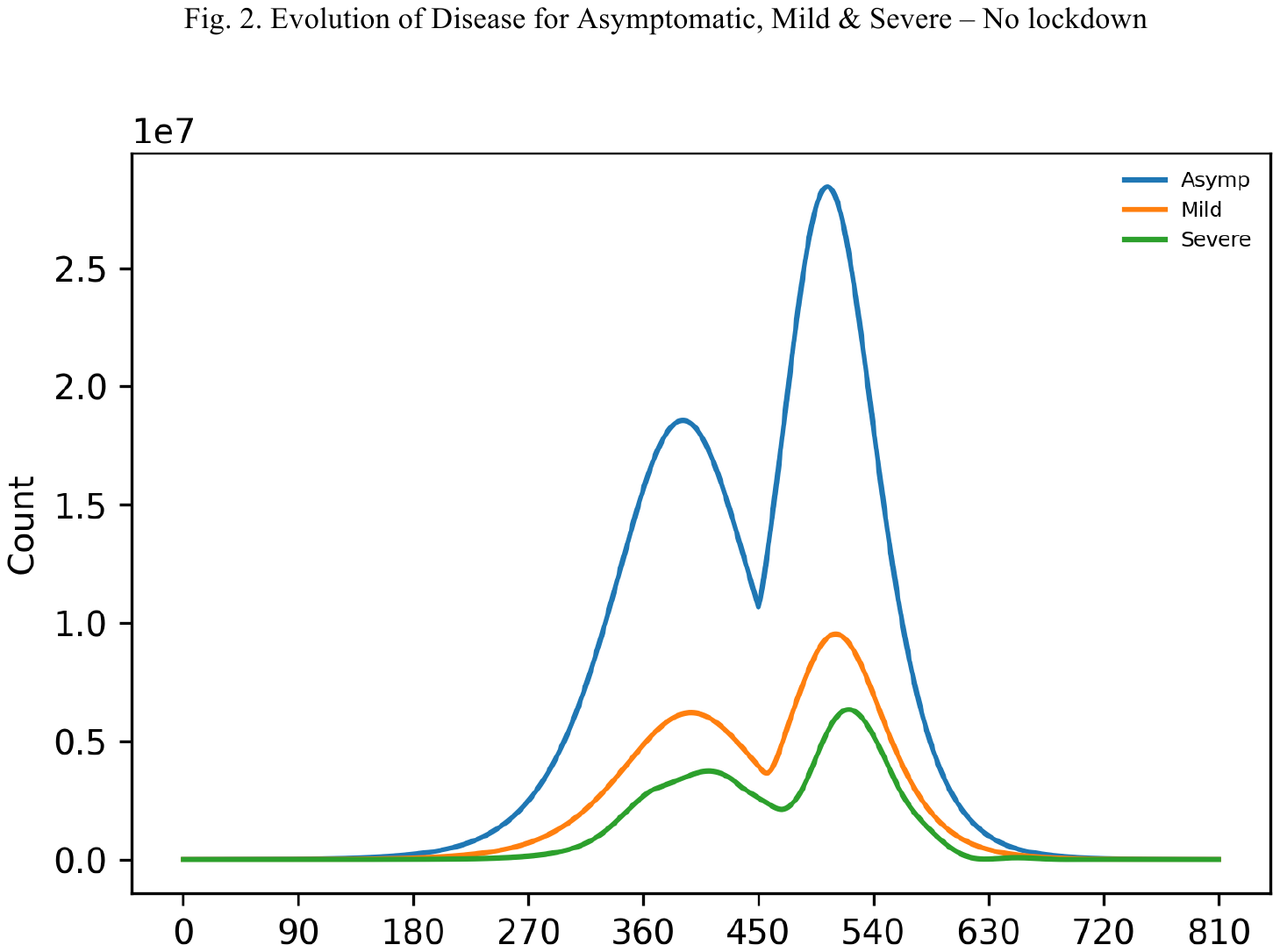}
\label{fig2}
\end{figure}

\begin{figure}[H]
\vskip-2cm
\hskip-2cm
\includegraphics[width=1.4\textwidth]{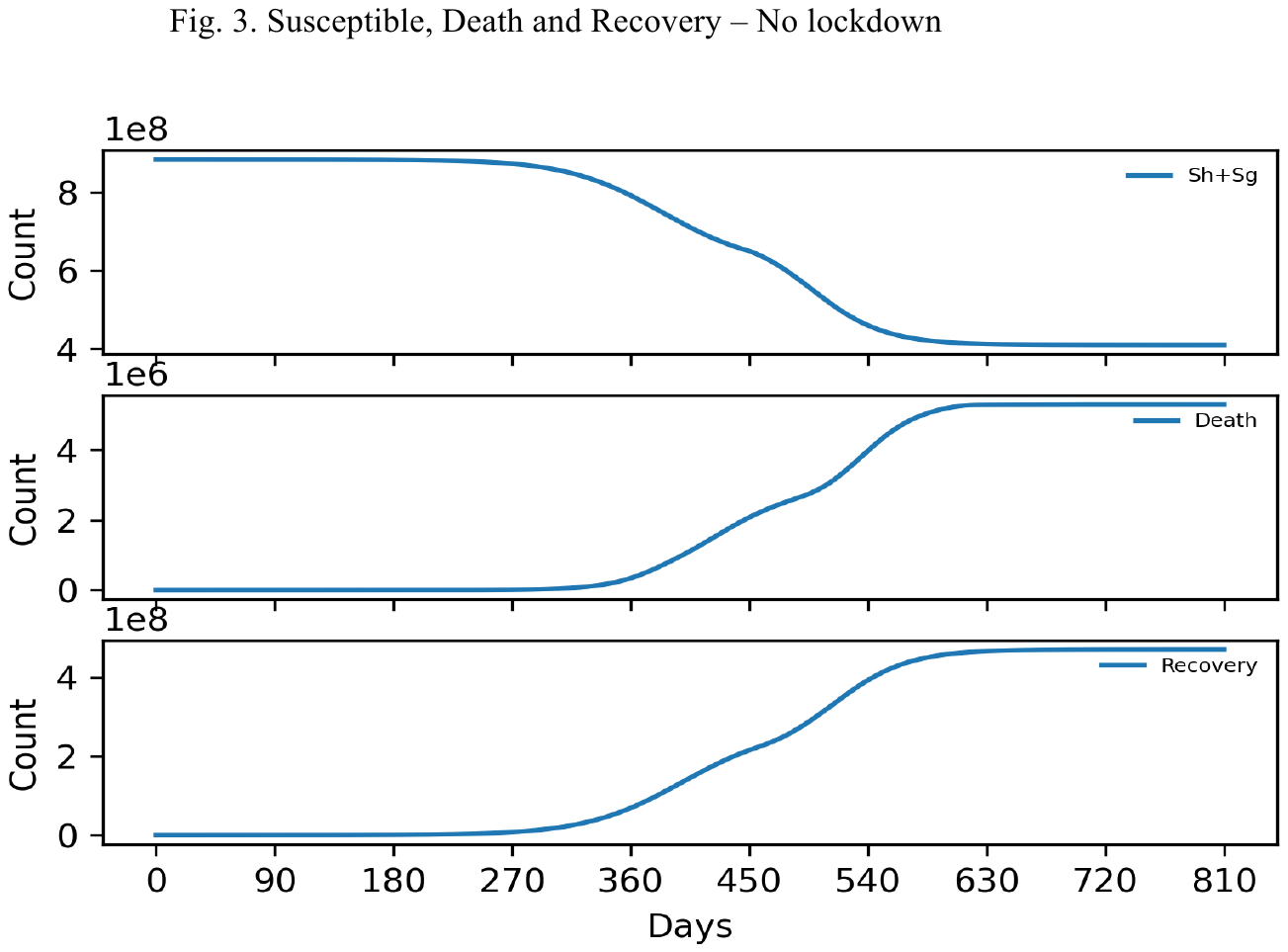}
\label{fig3}
\end{figure}

\begin{figure}[H]
\vskip-2cm
\hskip-2cm
\includegraphics[width=1.4\textwidth]{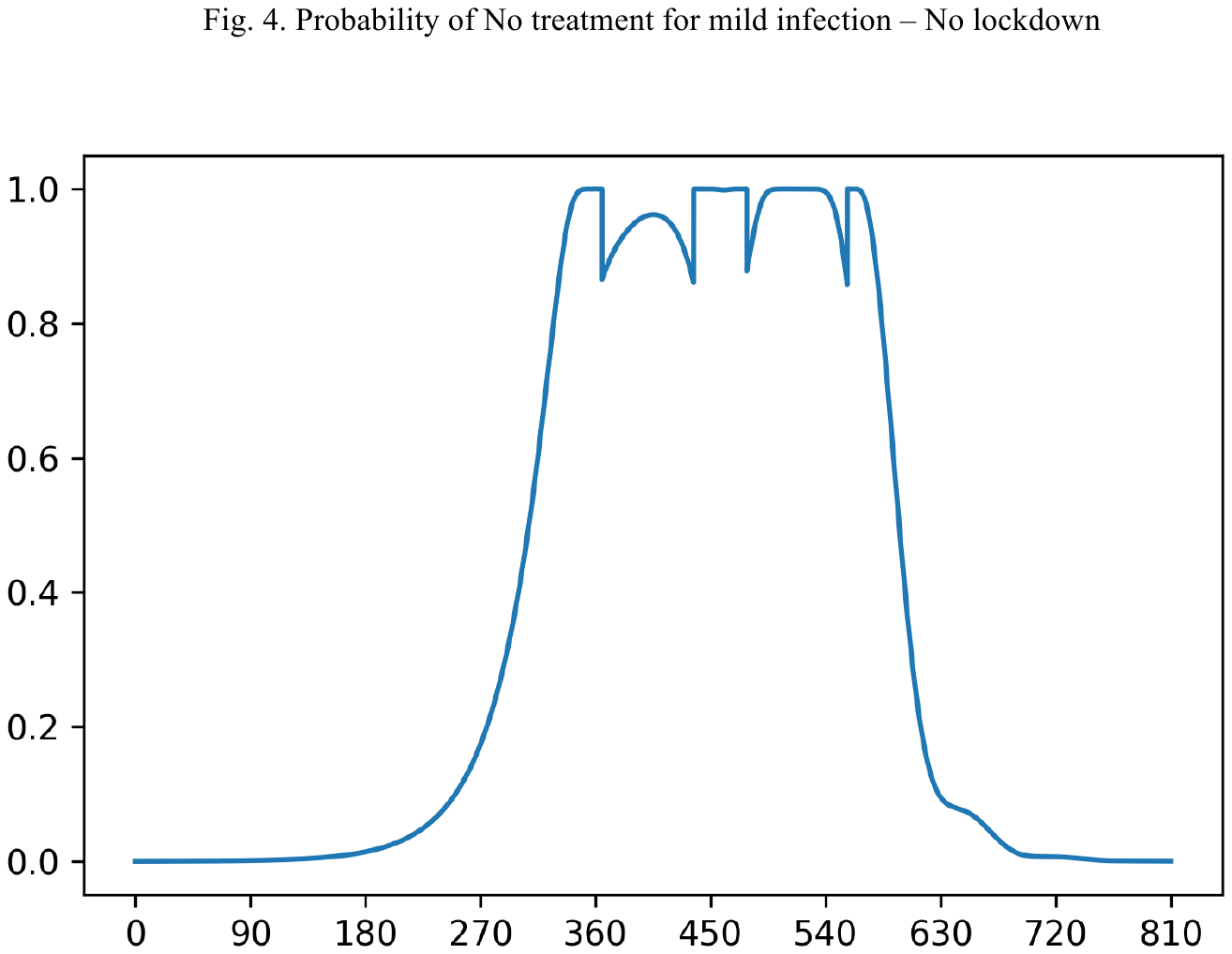}
\label{fig4}
\end{figure}

\begin{figure}[H]
\vskip-2cm
\hskip-2cm
\includegraphics[width=1.4\textwidth]{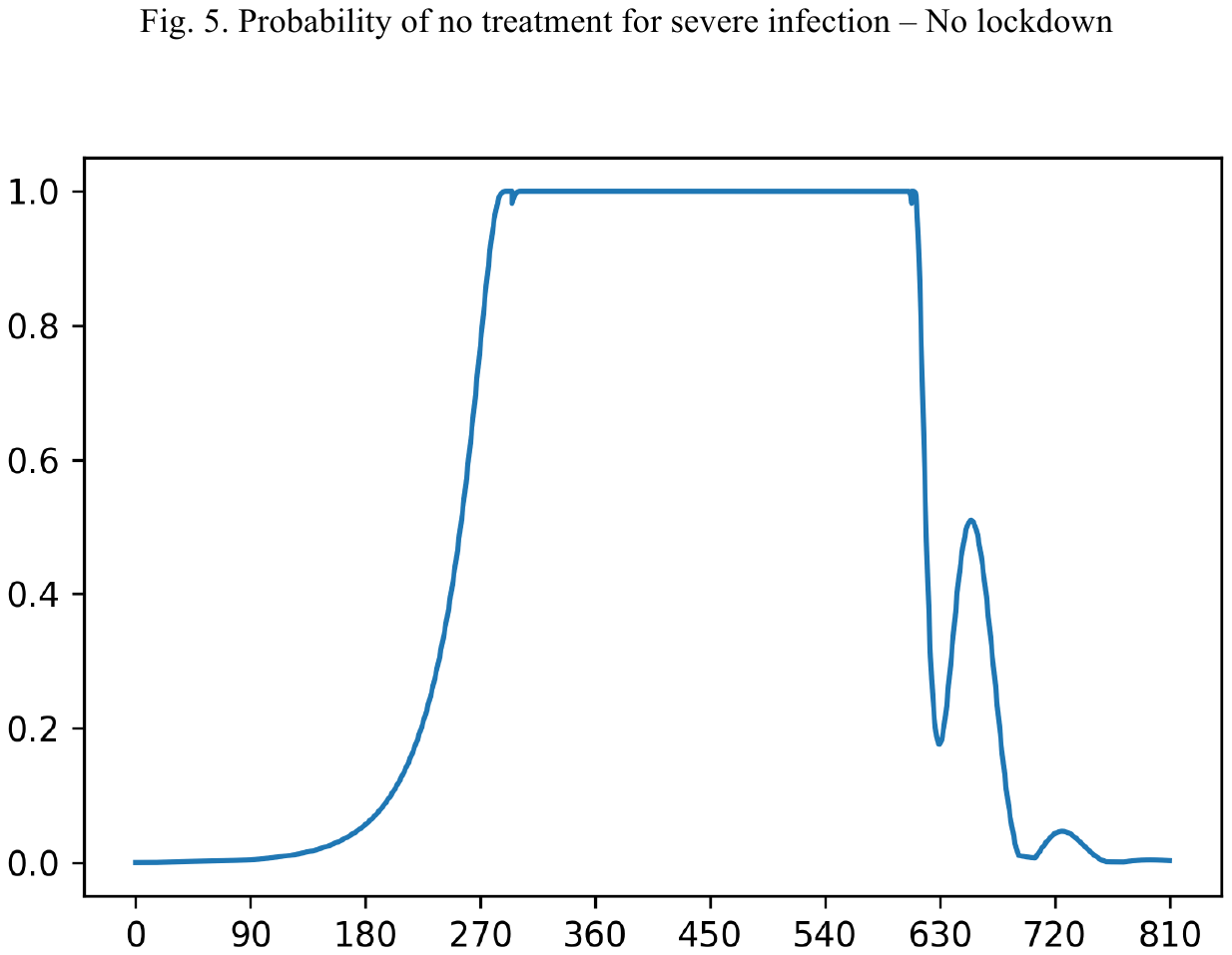}
\label{fig5}
\end{figure}

\begin{figure}[H]
\vskip-2cm
\hskip-2cm
\includegraphics[width=1.4\textwidth]{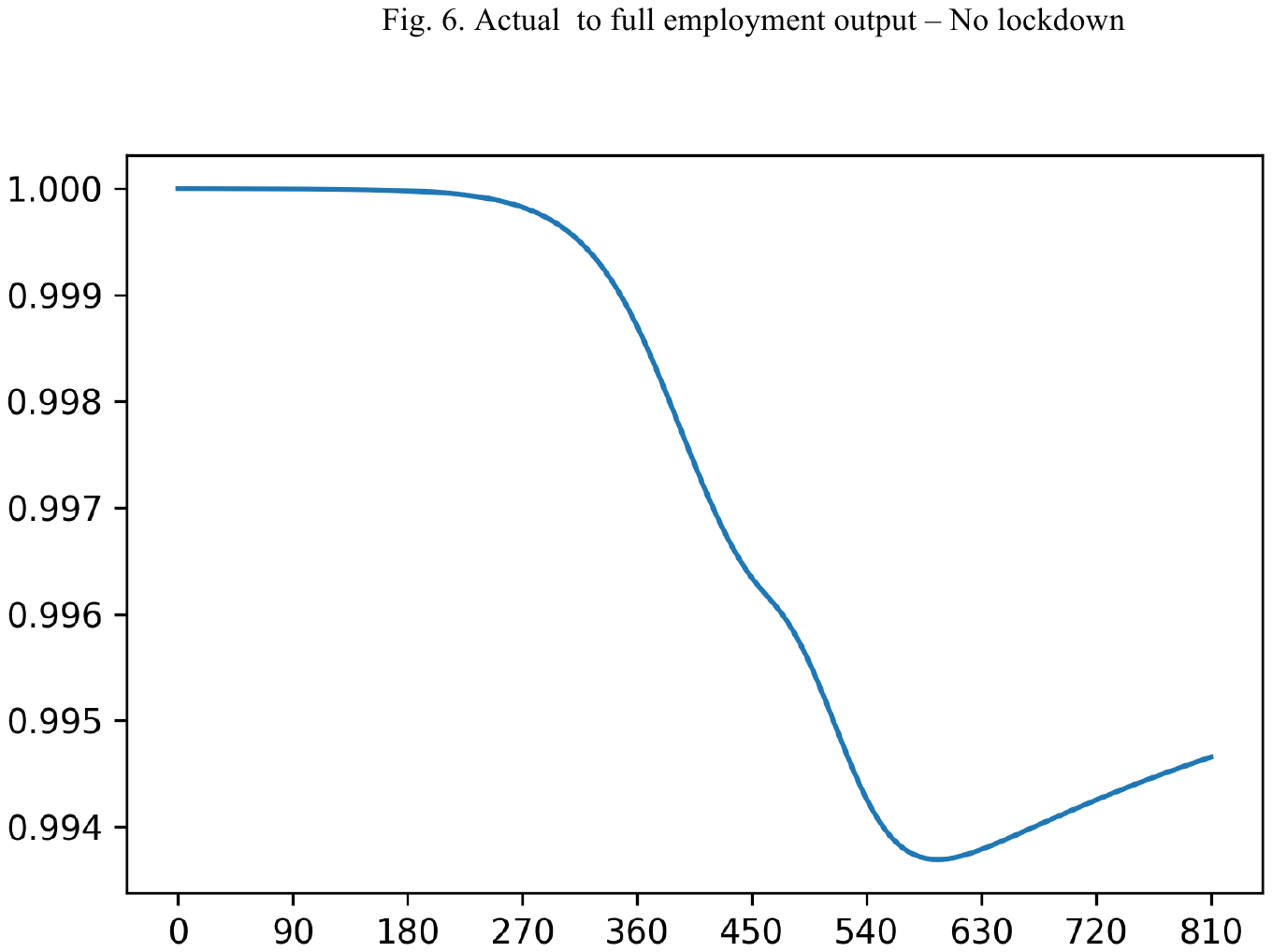}
\label{fig6}
\end{figure}

\begin{figure}[H]
\vskip-2cm
\hskip-2cm
\includegraphics[width=1.4\textwidth]{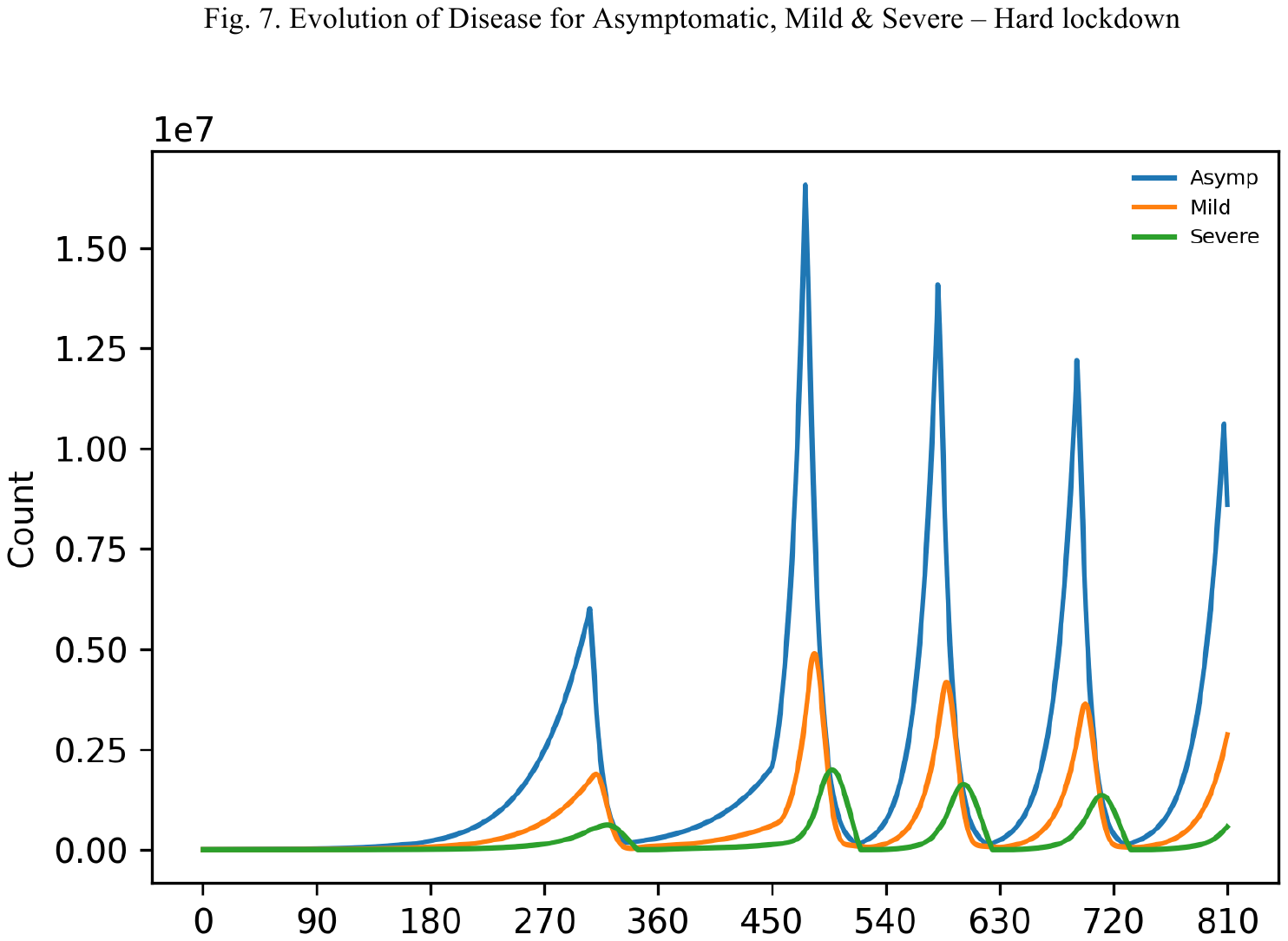}
\label{fig7}
\end{figure}

\begin{figure}[H]
\vskip-2cm
\hskip-2cm
\includegraphics[width=1.4\textwidth]{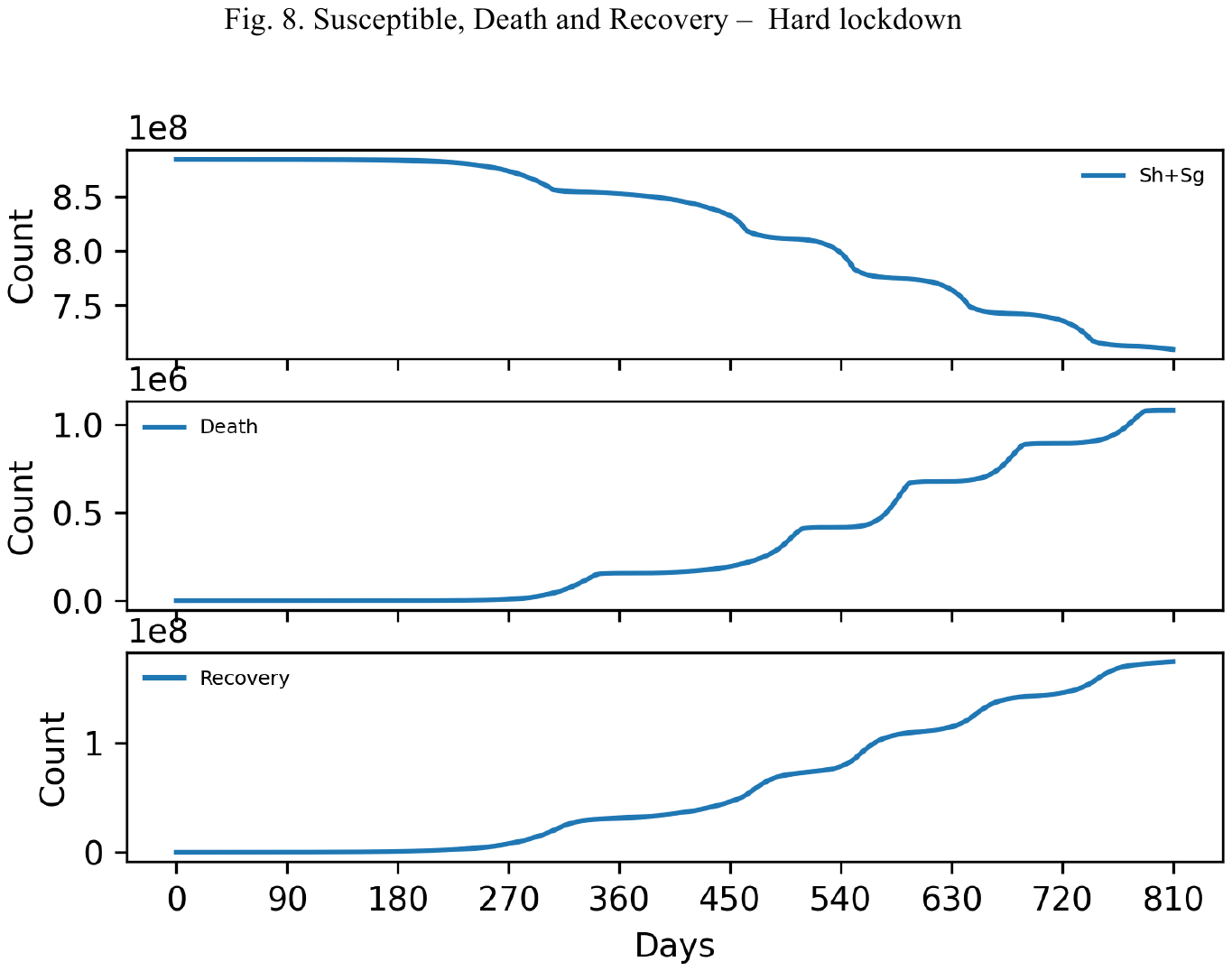}
\label{fig8}
\end{figure}

\begin{figure}[H]
\vskip-2cm
\hskip-2cm
\includegraphics[width=1.4\textwidth]{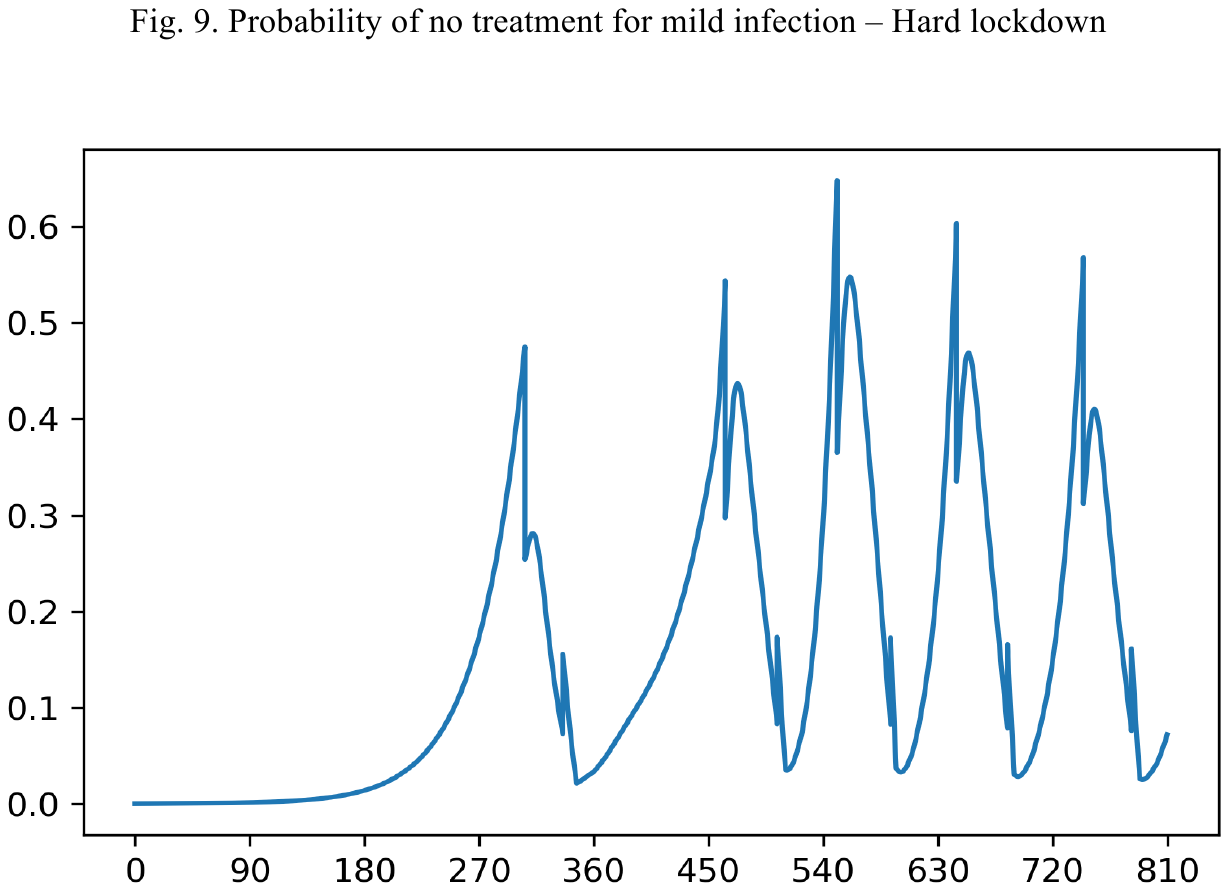}
\label{fig9}
\end{figure}

\begin{figure}[H]
\vskip-2cm
\hskip-2cm
\includegraphics[width=1.4\textwidth]{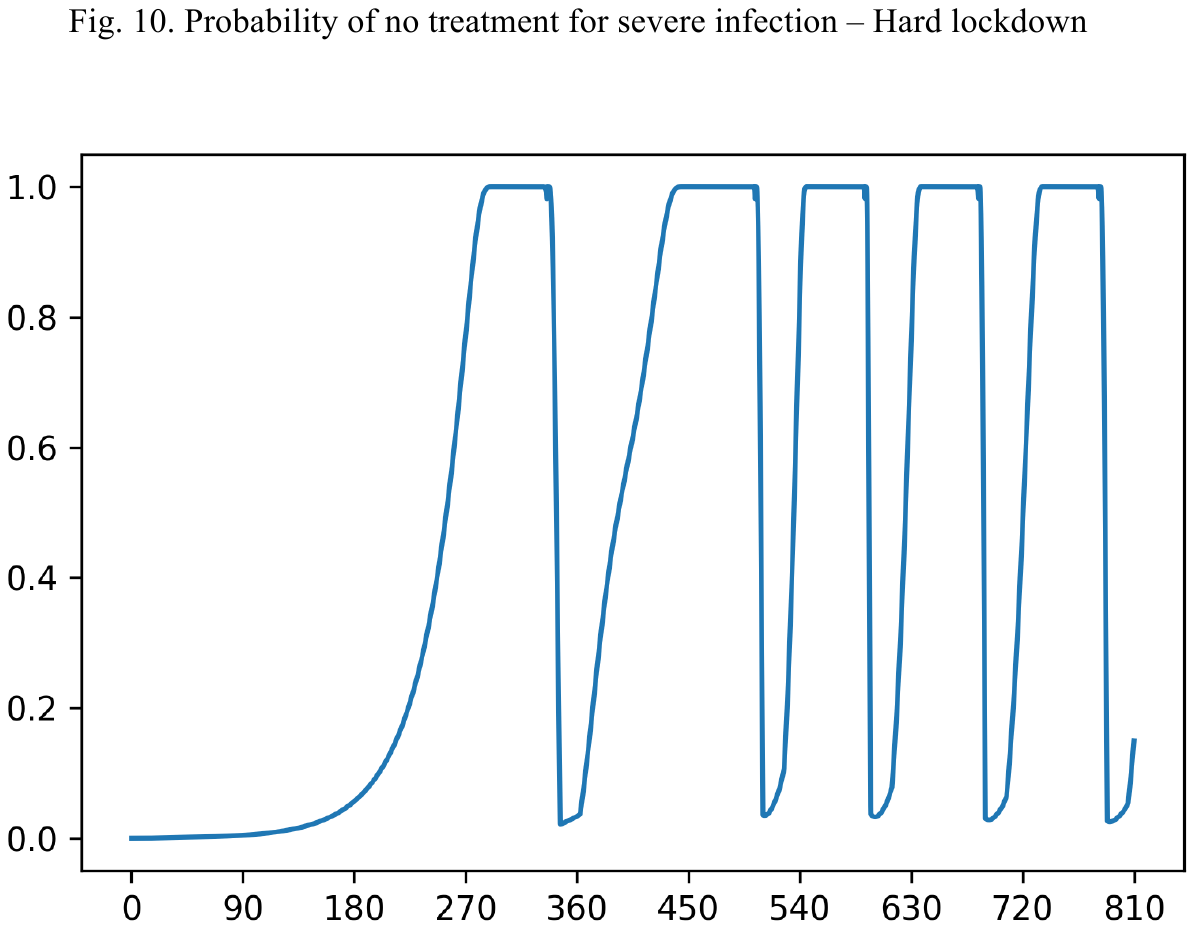}
\label{fig10}
\end{figure}

\begin{figure}[H]
\vskip-2cm
\hskip-2cm
\includegraphics[width=1.4\textwidth]{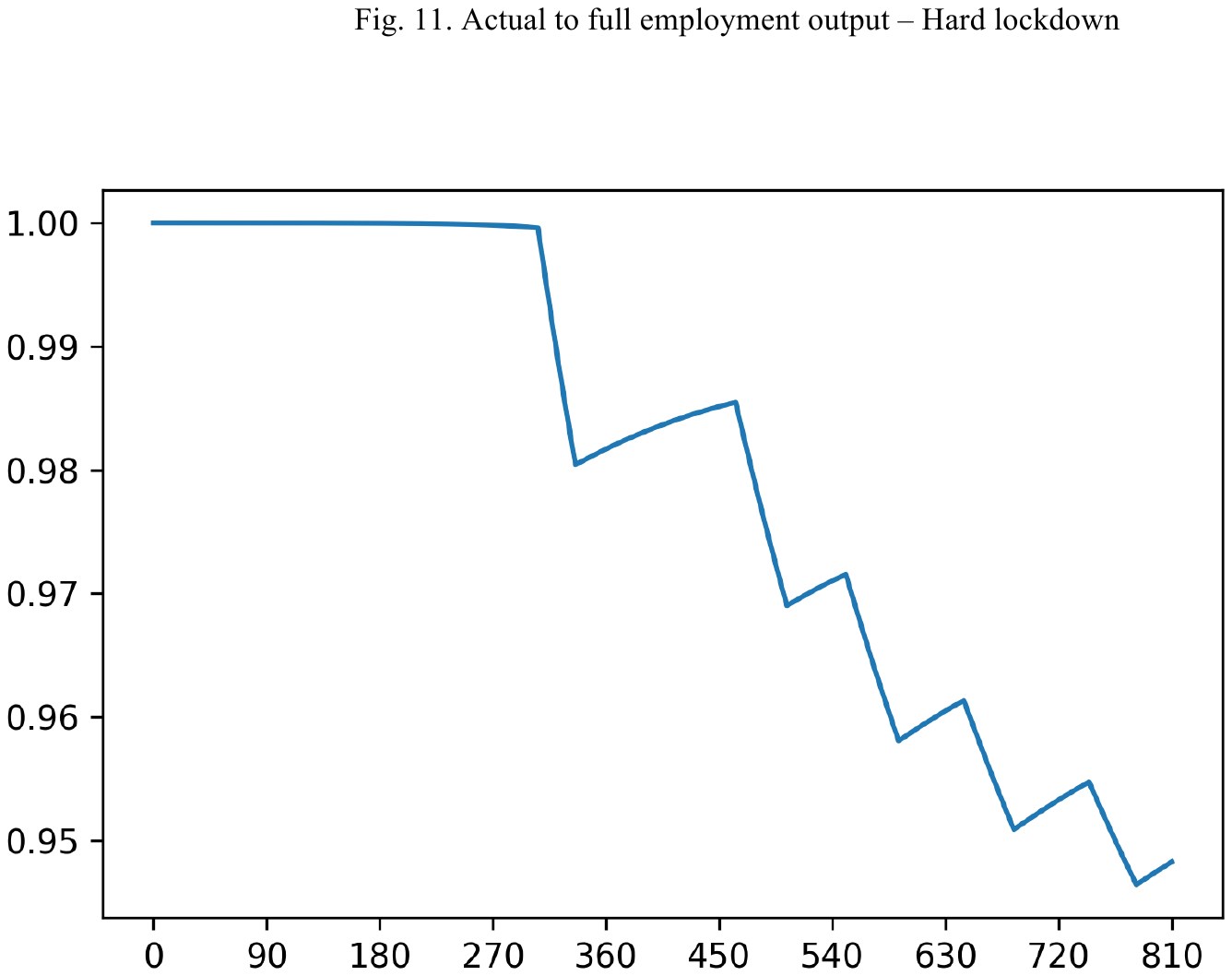}
\label{fig11}
\end{figure}

\begin{figure}[H]
\vskip-2cm
\hskip-2cm
\includegraphics[width=1.4\textwidth]{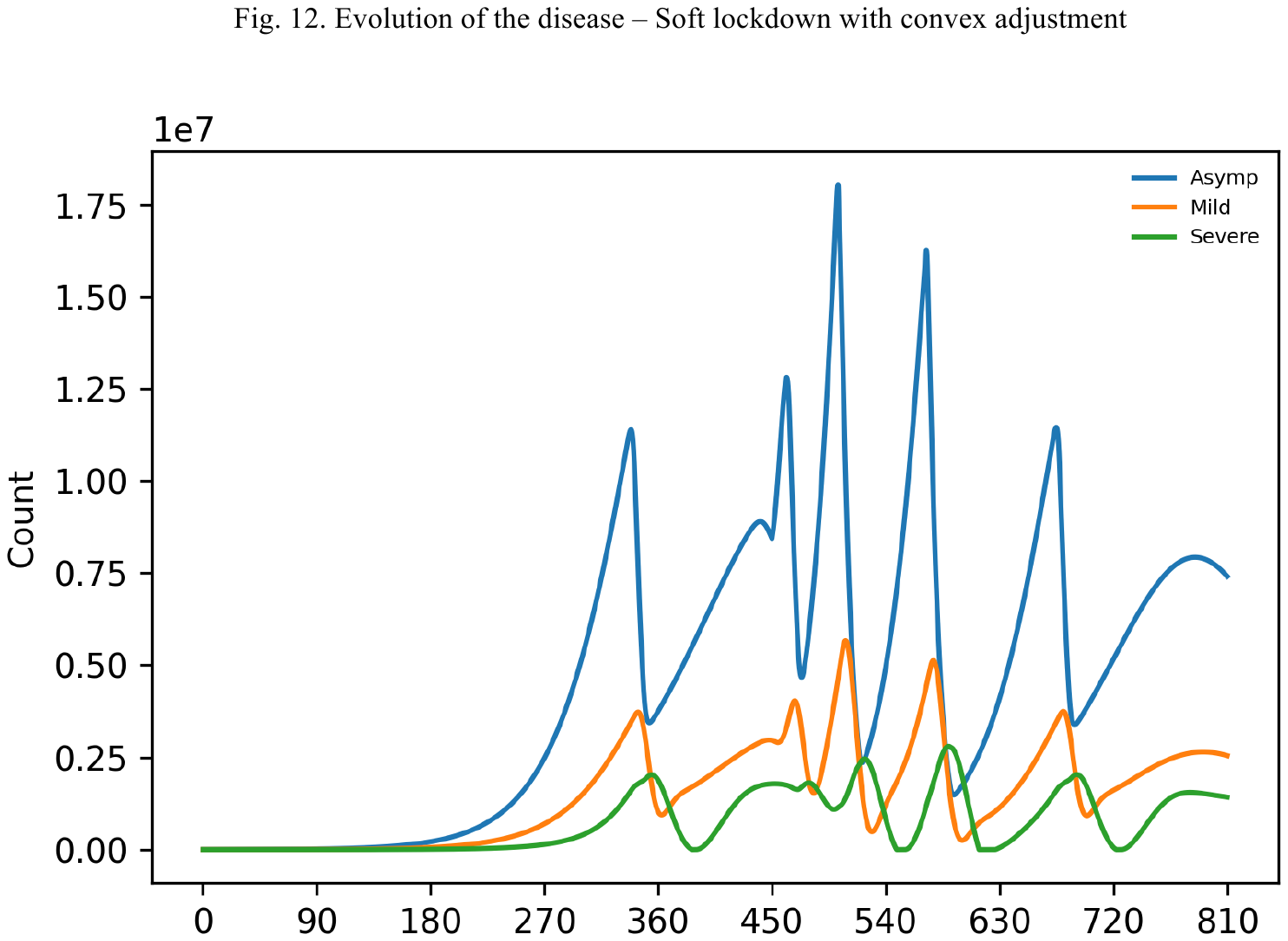}
\label{fig12}
\end{figure}

\begin{figure}[H]
\vskip-2cm
\hskip-2cm
\includegraphics[width=1.4\textwidth]{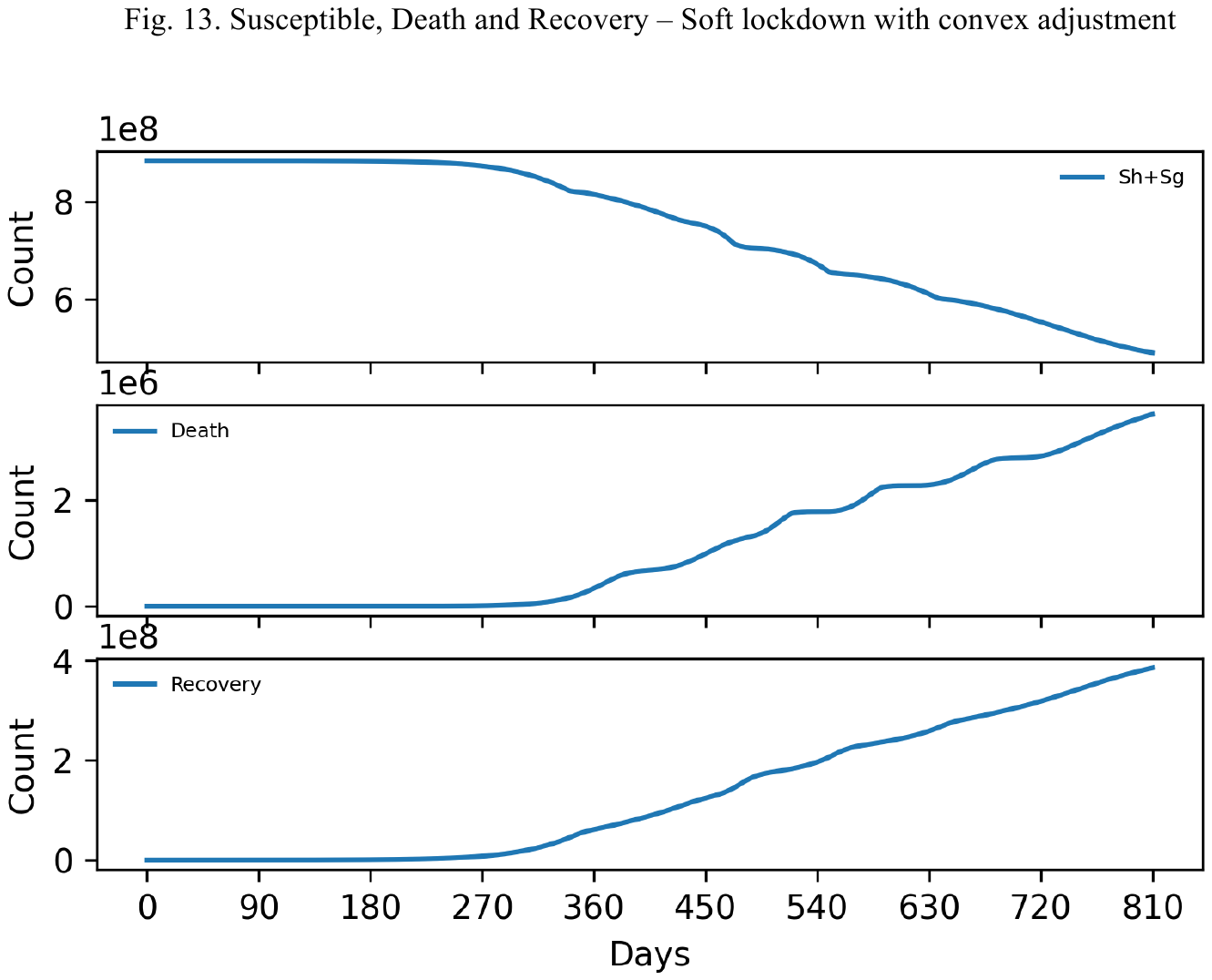}
\label{fig13}
\end{figure}

\begin{figure}[H]
\vskip-5cm
\hskip-2cm
\includegraphics[width=1.4\textwidth]{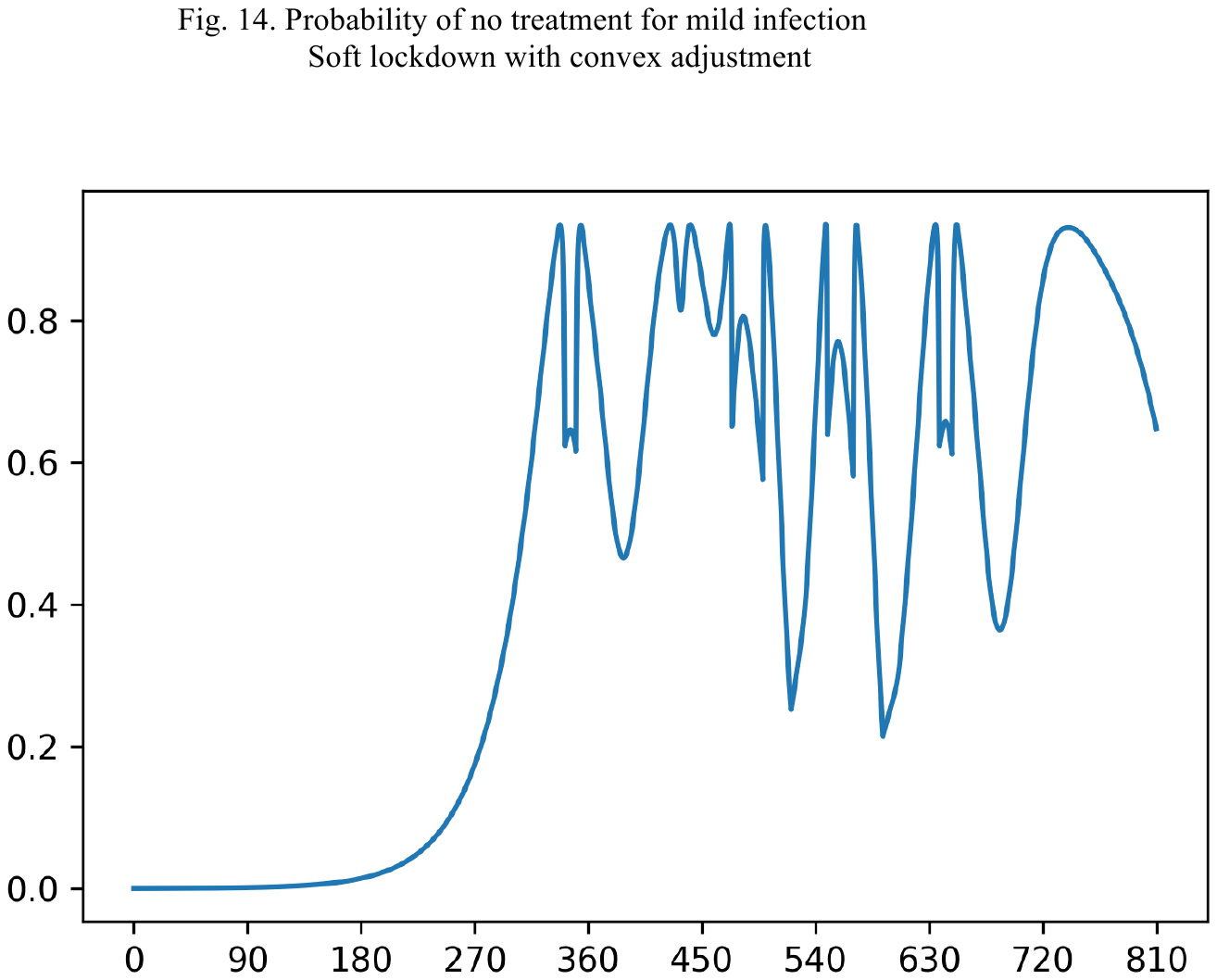}
\label{fig14}
\end{figure}

\begin{figure}[H]
\vskip-2cm
\hskip-2cm
\includegraphics[width=1.4\textwidth]{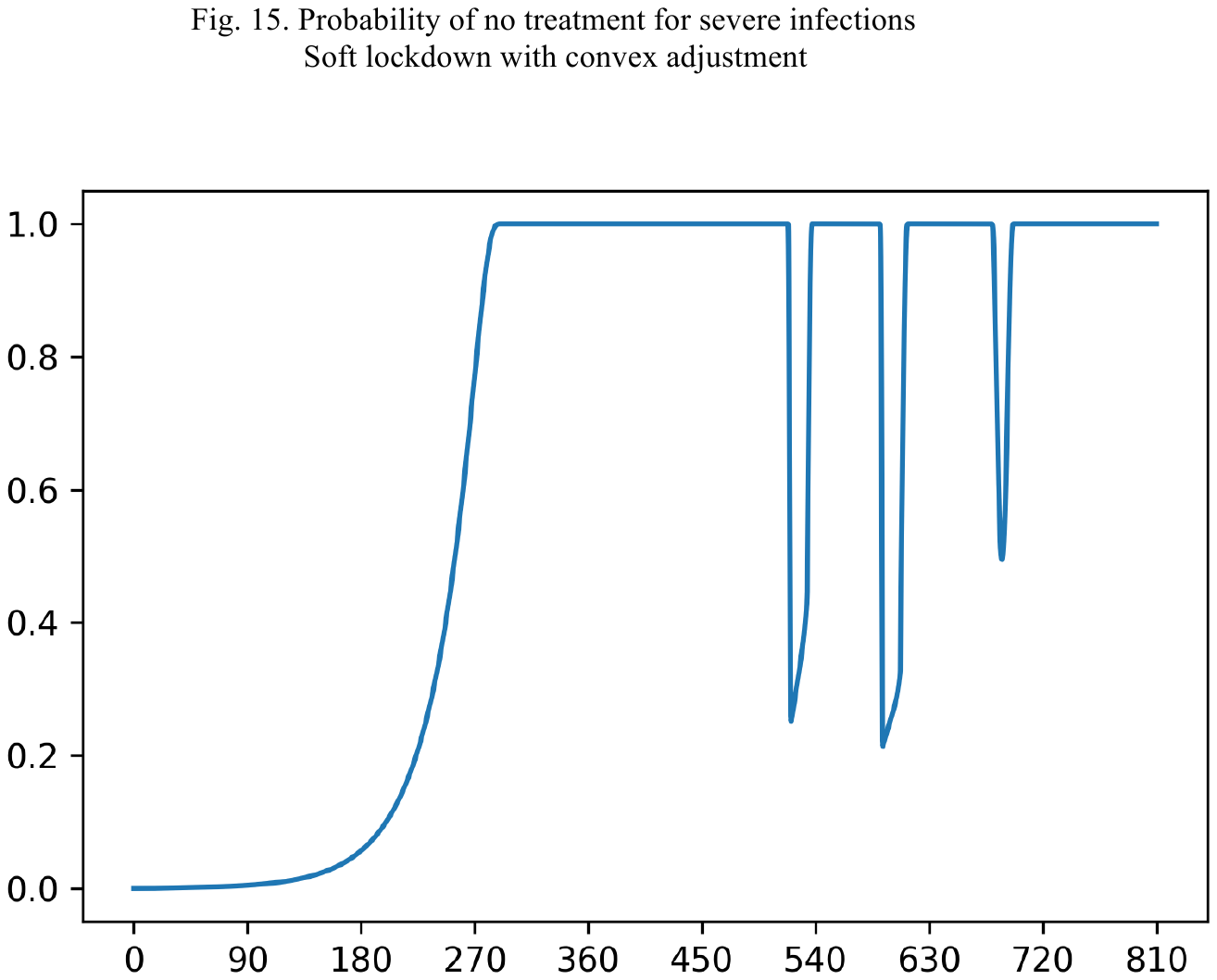}
\label{fig15}
\end{figure}

\begin{figure}[H]
\vskip-2cm
\hskip-2cm
\includegraphics[width=1.4\textwidth]{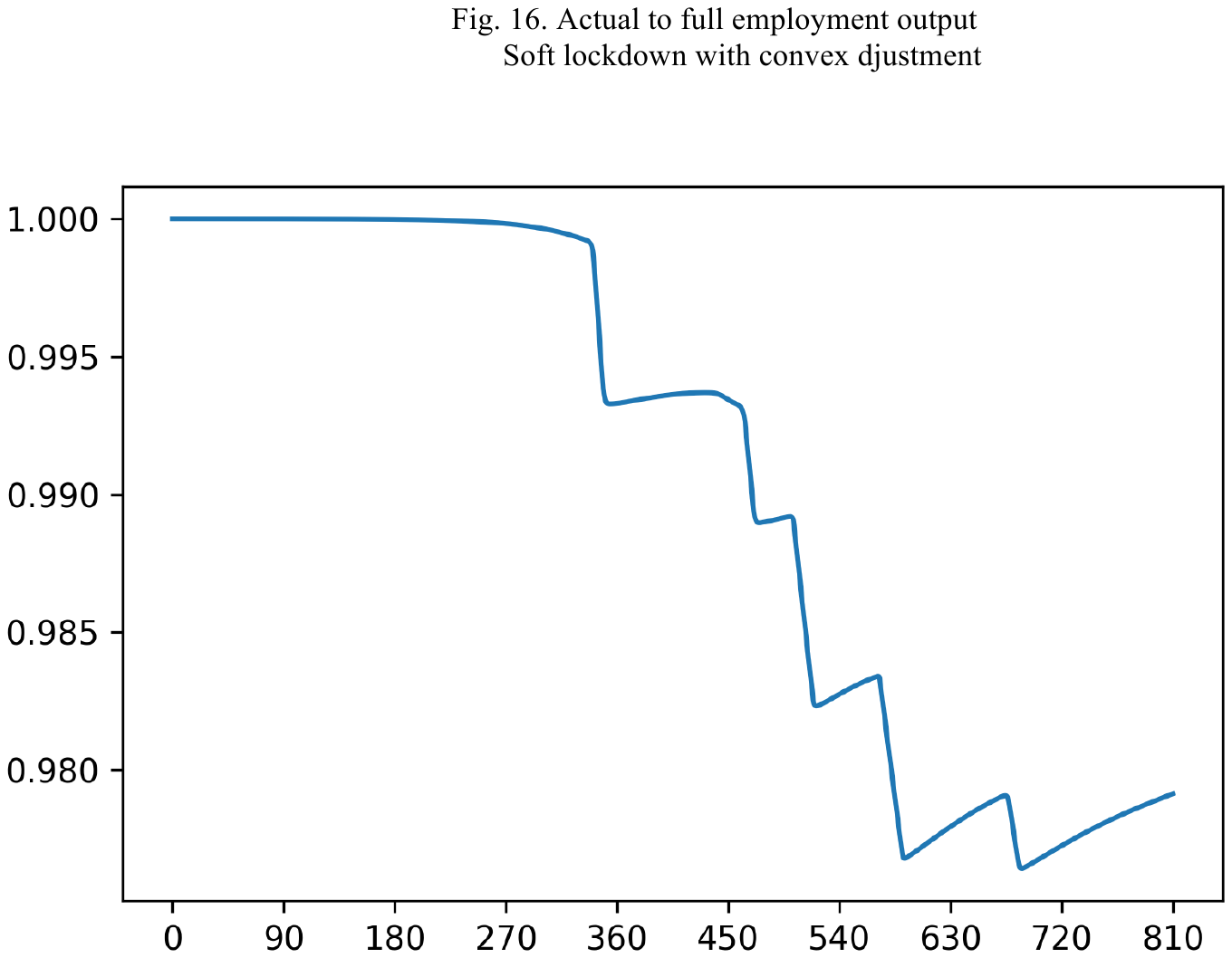}
\label{fig16}
\end{figure}

\begin{figure}[H]
\vskip-2cm
\hskip-2cm
\includegraphics[width=1.4\textwidth]{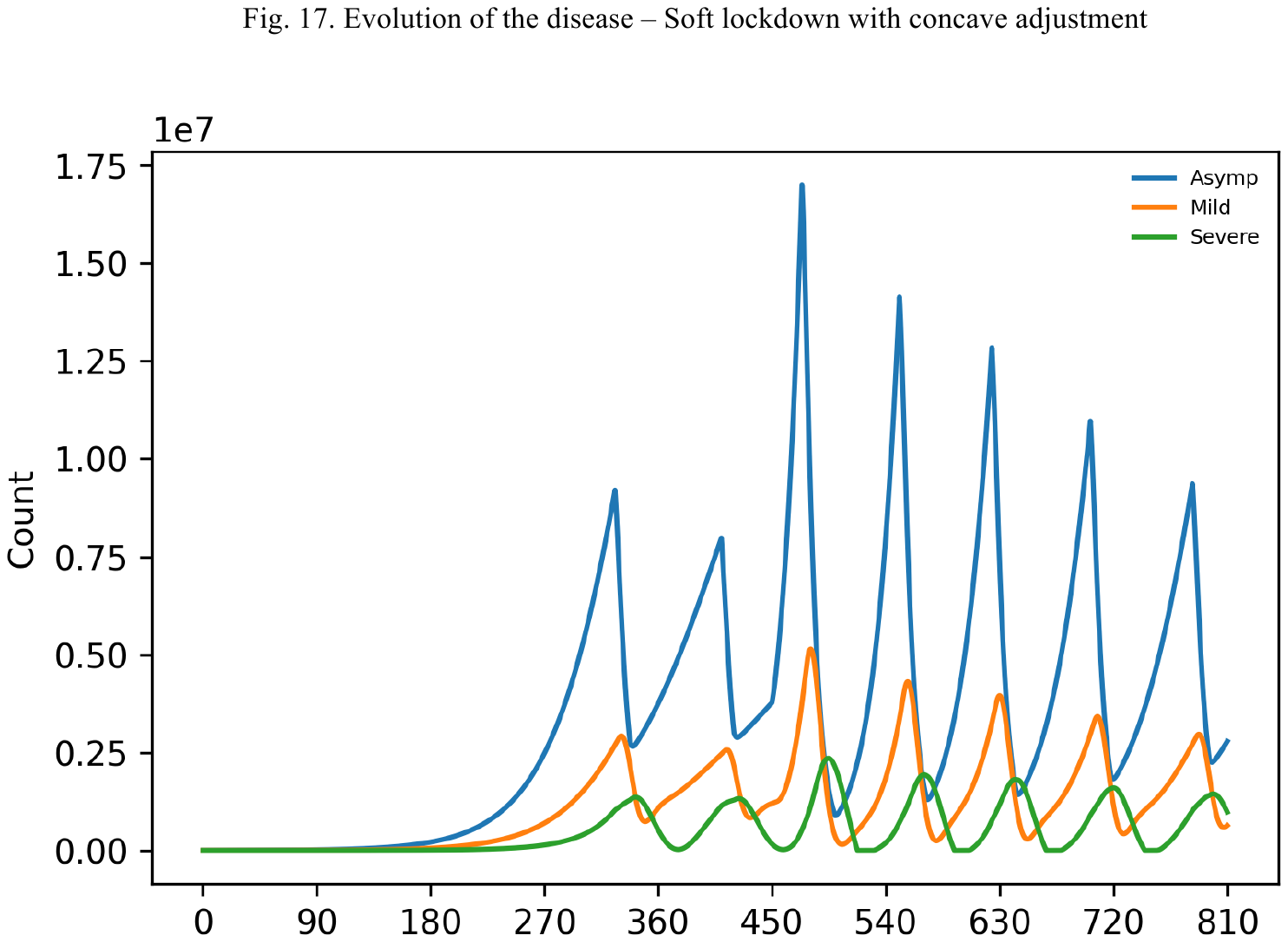}
\label{fig17}
\end{figure}

\begin{figure}[H]
\vskip-2cm
\hskip-2cm
\includegraphics[width=1.4\textwidth]{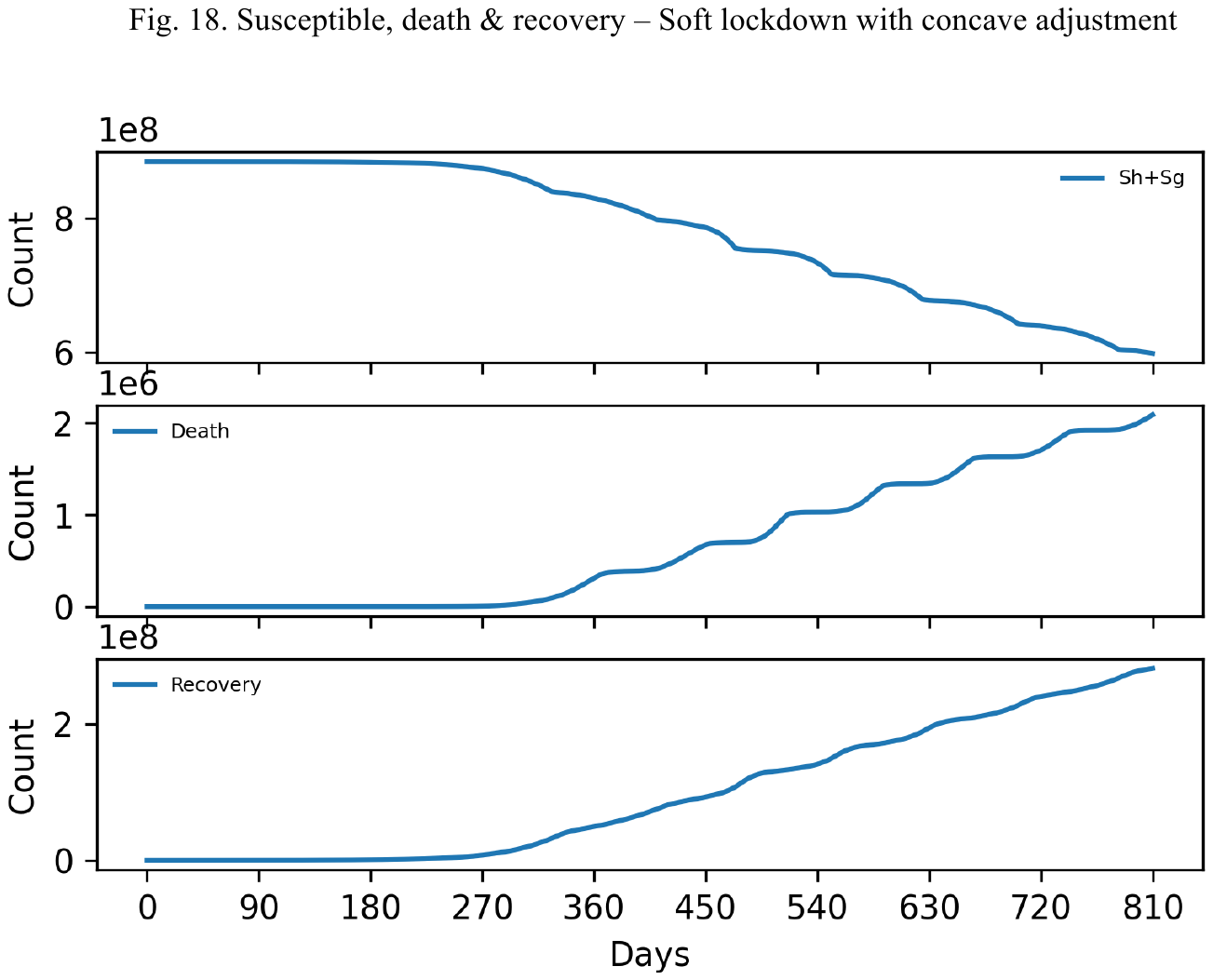}
\label{fig18}
\end{figure}

\begin{figure}[H]
\vskip-2cm
\hskip-2cm
\includegraphics[width=1.4\textwidth]{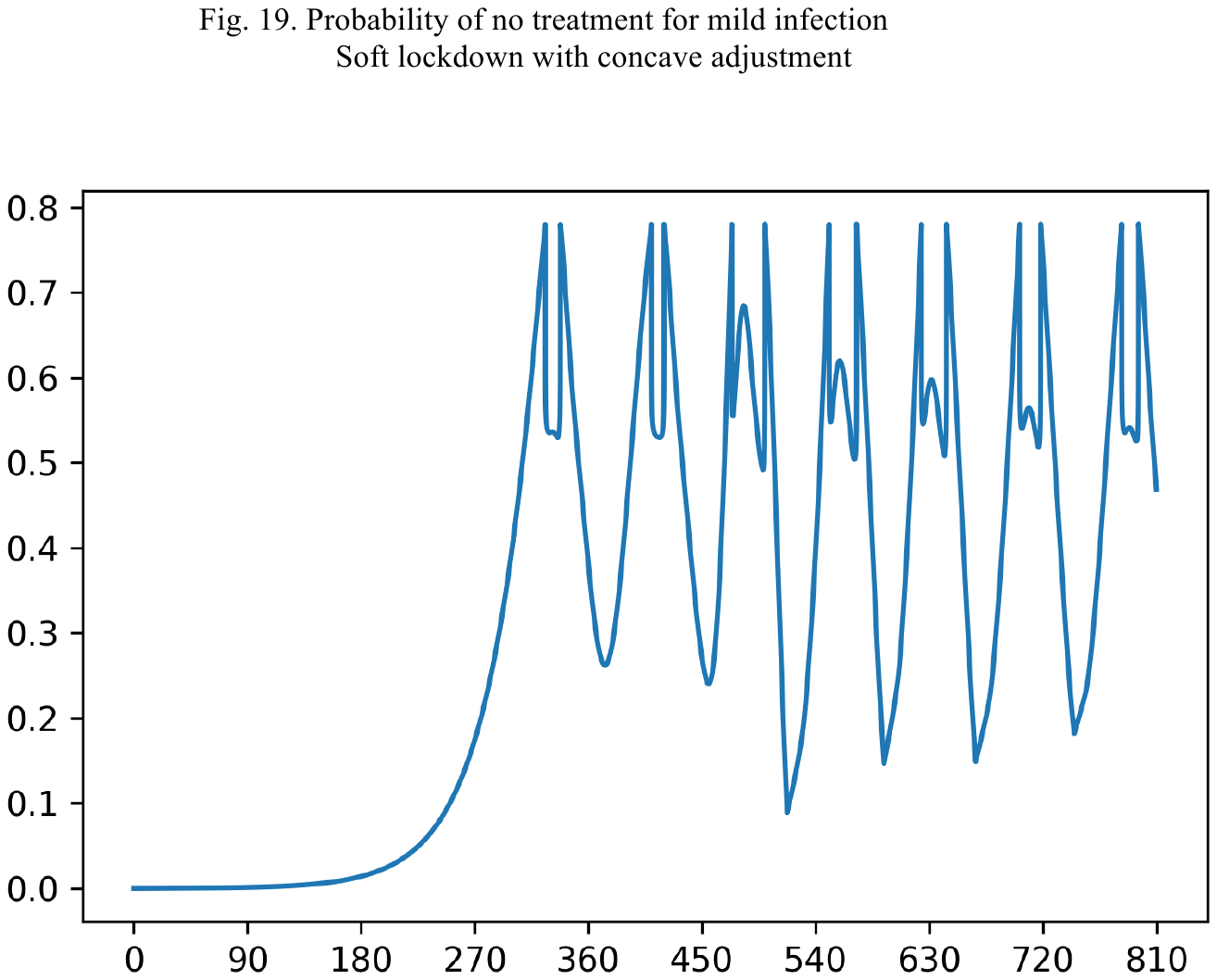}
\label{fig19}
\end{figure}

\begin{figure}[H]
\vskip-2cm
\hskip-2cm
\includegraphics[width=1.4\textwidth]{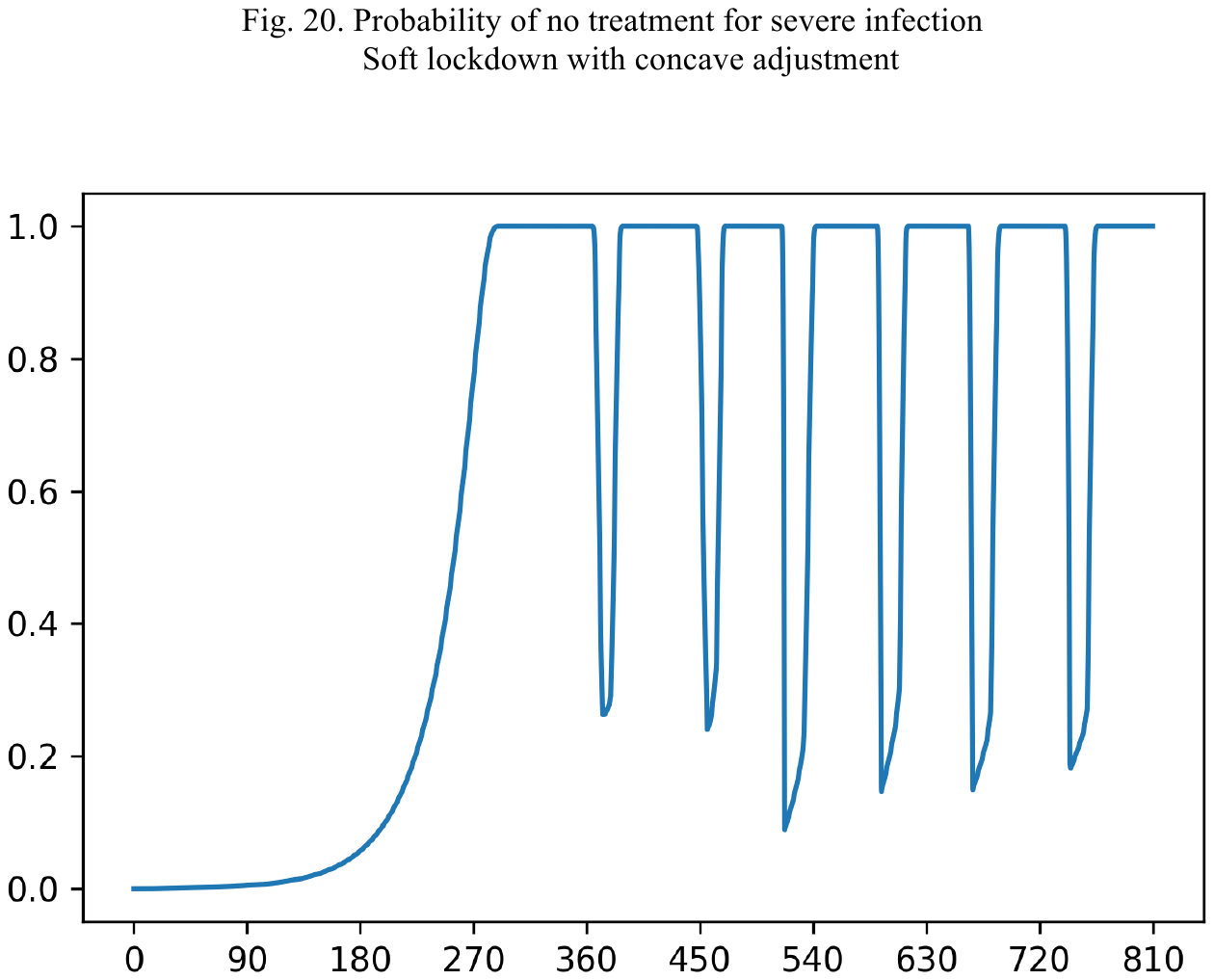}
\label{fig20}
\end{figure}

\begin{figure}[H]
\vskip-2cm
\hskip-2cm
\includegraphics[width=1.4\textwidth]{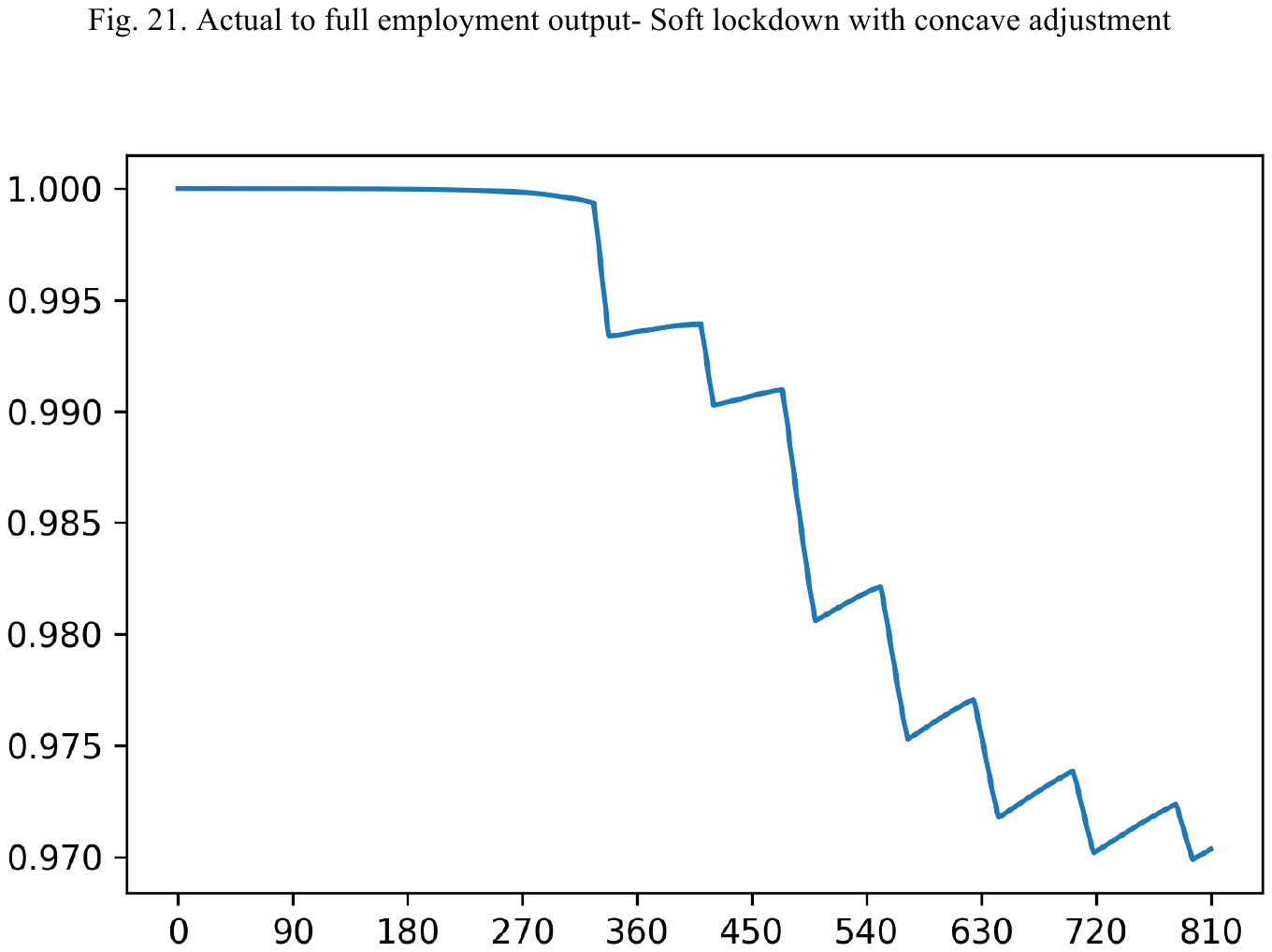}
\label{fig21}
\end{figure}

\begin{figure}[H]
\vskip-2cm
\hskip-2cm
\includegraphics[width=1.4\textwidth]{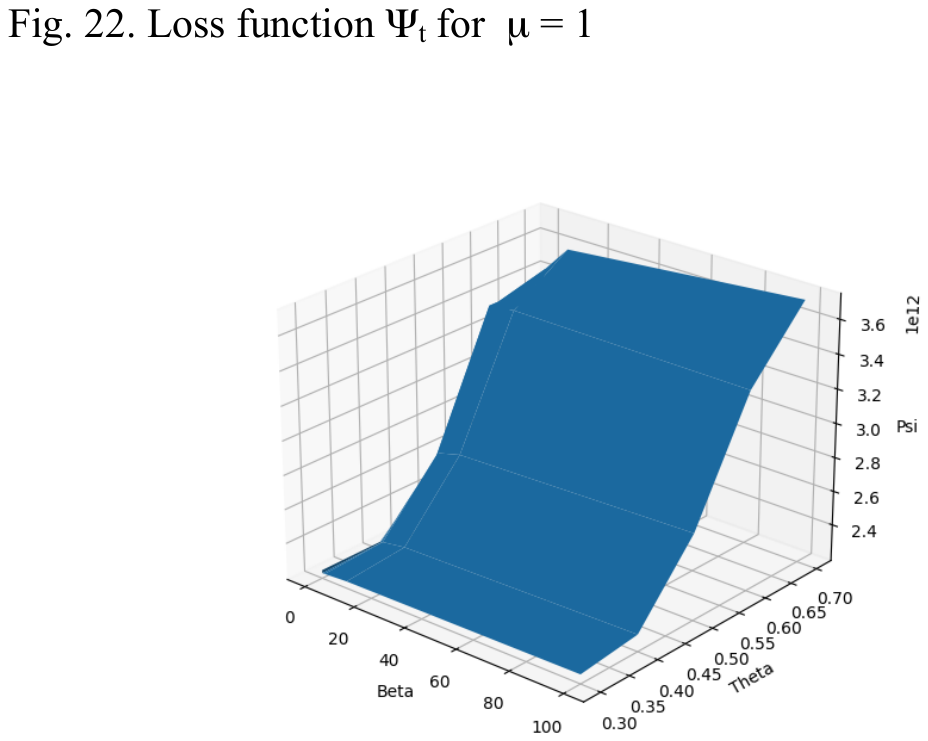}
\label{fig22}
\end{figure}

\newpage
\section*{Tables}

\begin{table}[H]
\centering
\begin{tiny}
\caption{\small Parameter values}
\begin{tabular}{cccccccccccccc}
\\
\hline
\multicolumn{1}{|c|}{$S_g$}       & \multicolumn{1}{c|}{$S_h$}       & \multicolumn{1}{c|}{$A$}                           & \multicolumn{1}{c|}{$\alpha_h$}                    & \multicolumn{1}{c|}{$B$}                           & \multicolumn{1}{c|}{$\beta_0$} 
&
\multicolumn{1}{c|}{$\alpha_0$}  & \multicolumn{1}{c|}{$\beta_1$}   & \multicolumn{1}{c|}{$\alpha_1$}                    & \multicolumn{1}{c|}{$\alpha_{22}$}                 & \multicolumn{1}{c|}{$\alpha_{42}$}                 & \multicolumn{1}{c|}{$\gamma_1$} & \multicolumn{1}{c|}{$\alpha_3$} 
& \multicolumn{1}{c|}{} 
\\ \hline
\multicolumn{1}{|c|}{884e6}   & \multicolumn{1}{c|}{0.76e6}      & \multicolumn{1}{c|}{1000}                          & \multicolumn{1}{c|}{3}                             & \multicolumn{1}{c|}{0.49e6}                        & \multicolumn{1}{c|}{0.01}  
& \multicolumn{1}{c|}{0.4}         & \multicolumn{1}{c|}{0.13}        & \multicolumn{1}{c|}{0.1}                           & \multicolumn{1}{c|}{0.3}                           & \multicolumn{1}{c|}{0.1}                           & \multicolumn{1}{c|}{0.06}       & \multicolumn{1}{c|}{0.01}
& \multicolumn{1}{c|}{} 
\\ \hline
                                                                                                                                                            &                                  &                                  &                                                    &                                                    &                                                    &                                 &                                 \\ \hline
\multicolumn{1}{|c|}{$\lambda_m$} & \multicolumn{1}{c|}{$\lambda_c$} & \multicolumn{1}{c|}{$\lambda_b$}                   & \multicolumn{1}{c|}{$\alpha$}                      & \multicolumn{1}{c|}{$c_m$}                         & \multicolumn{1}{c|}{$c_i$}      & \multicolumn{1}{c|}{$\nu$}  
& \multicolumn{1}{c|}{$\delta_h$}  & \multicolumn{1}{c|}{$\delta_g$}  & \multicolumn{1}{c|}{$
\delta_a
$} & \multicolumn{1}{c|}{$
\delta_m
$} & \multicolumn{1}{c|}
{$
\delta_i
$} & \multicolumn{1}{c|}{$\chi$}     & \multicolumn{1}{c|}{$\bar{Y}$}  \\ \hline
\multicolumn{1}{|c|}{2}           & \multicolumn{1}{c|}{2}           & \multicolumn{1}{c|}{2}                             & \multicolumn{1}{c|}{0.35}                          & \multicolumn{1}{c|}{0.02}                          & \multicolumn{1}{c|}{6}          & \multicolumn{1}{c|}{0.2}   
& \multicolumn{1}{c|}{14 days}     & \multicolumn{1}{c|}{14 days}     & \multicolumn{1}{c|}{3 days}                        & \multicolumn{1}{c|}{5 days}                        & \multicolumn{1}{c|}{14 days}                       & \multicolumn{1}{c|}{0.64e6}     & \multicolumn{1}{c|}{2.7e12}    
\\ \hline
                                                                              \end{tabular}
\end{tiny}
\label{tab1}
\end{table}

\begin{table}[H]
\centering
\begin{tiny}
\begin{tabular}{|c|c|c|}
\hline
Parameters  & No Lockdown & In Lockdown \\ \hline
$\lambda_0$ & 0.7         & 0.2         \\ \hline
$\lambda_1$ & 0.6         & 0.5         \\ \hline
\end{tabular}
\end{tiny}
\end{table}

\vskip2cm

\vskip-2cm
\hskip-3.0cm
\includegraphics[width=1.4\textwidth]{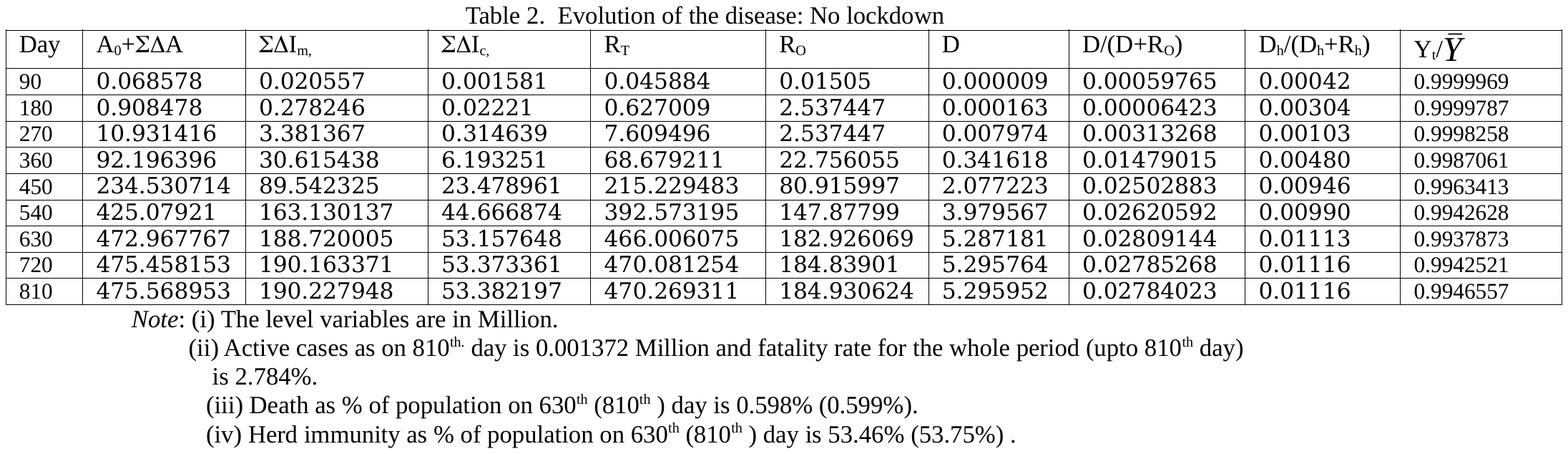}

\vskip-2cm
\hskip-3.2cm
\includegraphics[width=1.4\textwidth]{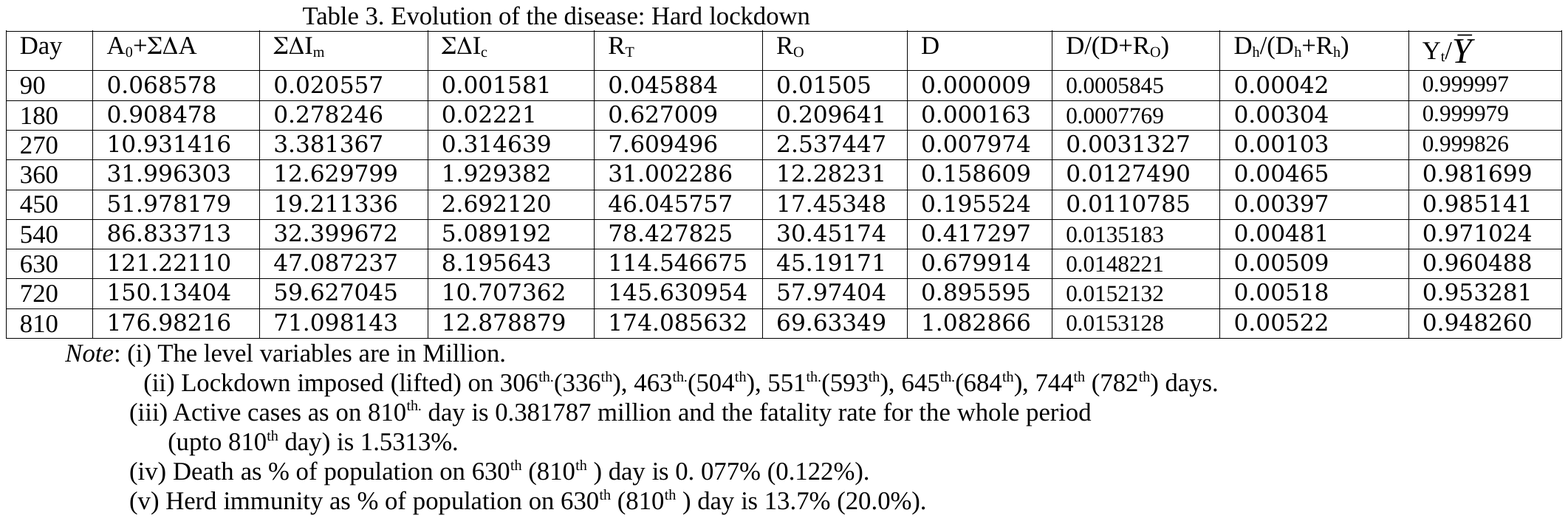}

\vskip-2.5cm
\hskip-3.8cm
\includegraphics[width=1.4\textwidth]{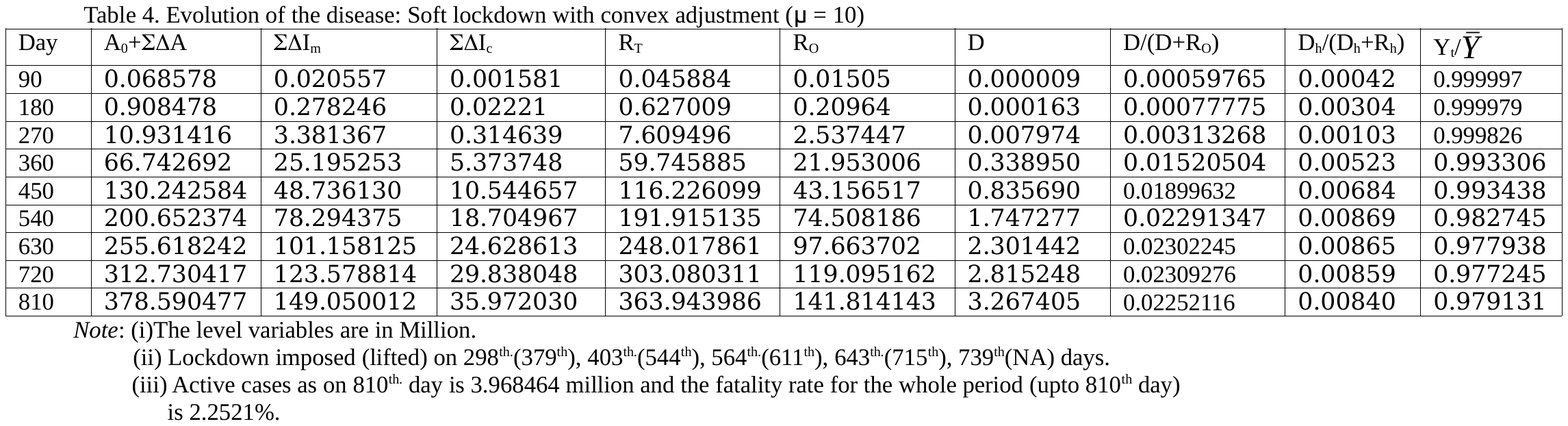}

\vskip-2.5cm
\hskip-3.7cm
\includegraphics[width=1.4\textwidth]{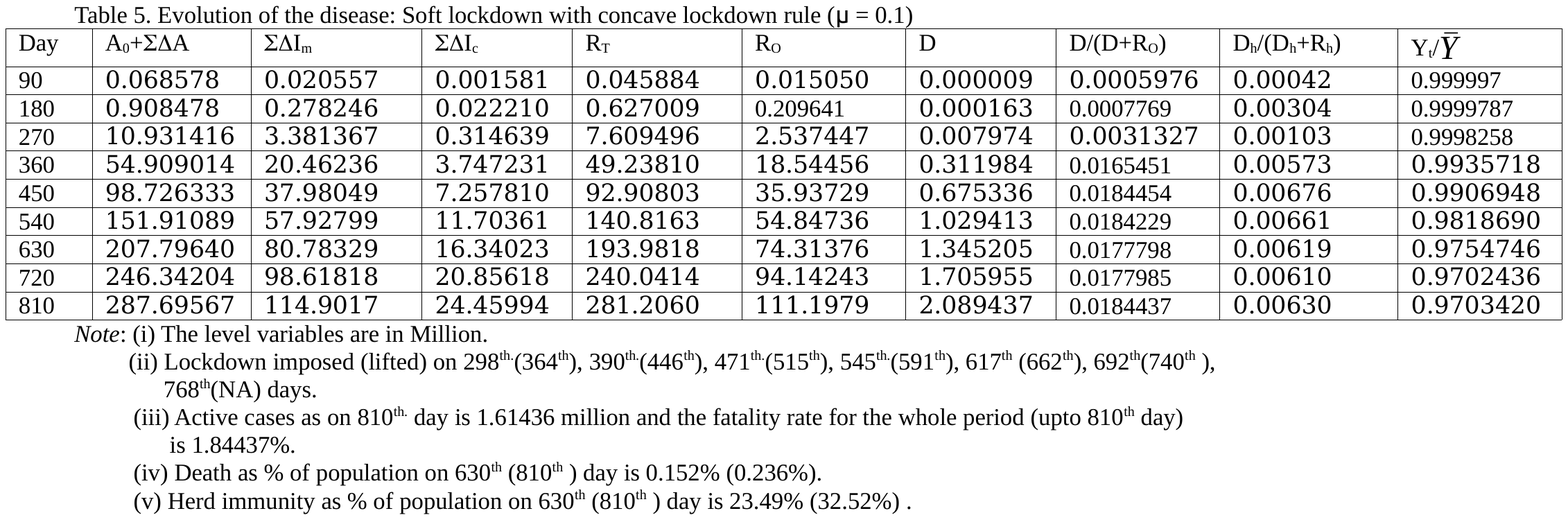}

\vskip-3cm
\hskip-2.2cm
\includegraphics[width=1.0\textwidth]{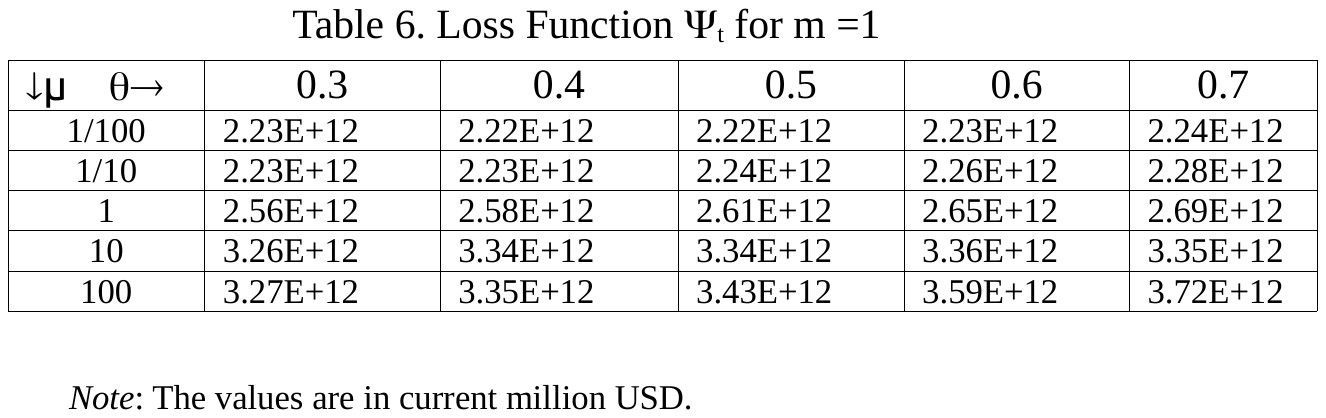}

\end{document}